\documentclass[namedate,webpdf,imanum]{ima-authoring-template}

\usepackage{bm}
\usepackage{amsfonts}
\usepackage{subcaption}
\usepackage[normalem]{ulem}

\newtheorem{lemma}{Lemma}
\newtheorem{corollary}{Corollary}
\newtheorem{remark}{Remark}

\DeclareMathOperator*{\argmax}{argmax}
\DeclareMathOperator*{\argmin}{argmin}
\DeclareMathOperator*{\esssup}{esssup}
\DeclareMathOperator*{\Real}{Re}
\DeclareMathOperator*{\diag}{diag}
\DeclareMathOperator*{\Span}{span}
\DeclareMathOperator{\supp}{supp}

\graphicspath{{Fig/}}

\theoremstyle{thmstyletwo}%
\newtheorem{theorem}{Theorem}%  meant for continuous numbers
\newtheorem{proposition}[theorem]{Proposition}%

\newtheorem{definition}{Definition}

\numberwithin{equation}{section}

\begin{document}

\DOI{DOI HERE}
\copyrightyear{2022}
\vol{00}
\pubyear{2021}
\access{Advance Access Publication Date: Day Month Year}
\appnotes{Paper}
\copyrightstatement{Published by Oxford University Press on behalf of the Institute of Mathematics and its Applications. All rights reserved.}
\firstpage{1}

%\subtitle{Subject Section}

\title[Compressive Fourier collocation methods]{Compressive Fourier collocation methods for high-dimensional diffusion equations with periodic boundary conditions}

\author{Weiqi Wang* \ORCID{0000-0002-2798-8181} and Simone Brugiapaglia \ORCID{0000-0003-1927-8232}
\address{\orgdiv{Department of Mathematics and Statistics}\\ \orgname{Concordia University}\\ \orgaddress{Montr\'eal, \state{QC}, \country{Canada}}}}

\authormark{Weiqi Wang and Simone Brugiapaglia}

\corresp[*]{Corresponding author: \href{email:weiqi.wang@concordia.ca}{weiqi.wang@concordia.ca}}

\received{Date}{0}{Year}
\revised{Date}{0}{Year}
\accepted{Date}{0}{Year}

%\editor{Associate Editor: Name}

\abstract{High-dimensional Partial Differential Equations (PDEs) are a popular mathematical modelling tool, with applications ranging from finance to computational chemistry. However, standard numerical techniques for solving these PDEs are typically affected by the curse of dimensionality. In this work, we tackle this challenge while focusing on stationary diffusion equations defined over a high-dimensional domain with periodic boundary conditions. Inspired by recent progress in sparse function approximation in high dimensions, we propose a new method called compressive Fourier collocation. Combining ideas from compressive sensing and spectral collocation, our method replaces the use of structured collocation grids with Monte Carlo sampling and employs sparse recovery techniques, such as orthogonal matching pursuit and $\ell^1$ minimization, to approximate the Fourier coefficients of the PDE solution. We conduct a rigorous theoretical analysis showing that the approximation error of the proposed method is comparable with the best $s$-term approximation (with respect to the Fourier basis) to the solution. Using the recently introduced framework of random sampling in bounded Riesz systems, our analysis shows that the compressive Fourier collocation method mitigates the curse of dimensionality with respect to the number of collocation points under sufficient conditions on the regularity of the diffusion coefficient. We also present numerical experiments that illustrate the accuracy and stability of the method for the approximation of sparse and compressible solutions.}

% \boxedtext{
% \begin{itemize}
% \item Key boxed text here.
% \item Key boxed text here.
% \item Key boxed text here.
% \end{itemize}}

\maketitle

\section{Introduction}

The numerical solution of high-dimensional Partial Differential Equations (PDEs) is a crucial task in scientific computing with a wide range of applications. Popular high-dimensional PDEs include the Black-Scholes equation in computational finance, the many-electron Schr\"odinger equation in computational chemistry, the Hamilton-Jacobi-Bellman equation in optimal control, and the Fokker-Planck equation in statistical mechanics. Computing numerical solutions to high-dimensional PDEs is made intrinsically challenging by the \emph{curse of dimensionality} (term coined by \citet{bellman1957dynamic,bellman1961adaptive}), i.e.\ the tendency of numerical techniques to require a computational cost that scales exponentially with respect to the dimension of the PDE domain. %Start for classical methods based on sparse grids, techniques base don low rank decompositions have been proposed, and more recently, the use of deep learning to numerically solve high-dimensional PDEs has attracted a huge amount of interest.  See \cite{weinan2021algorithms} for a recent review. 

We tackle this challenge by proposing a new numerical method for high-dimensional PDEs called \emph{compressive Fourier collocation}, based on ideas from compressive sensing. Our work is inspired by recent progress in sparse polynomial approximation of high-dimensional functions, where methods based on compressive sensing can provably lessen the curse of dimensionality with respect to the sample complexity (see \citet{adcock2022sparse} for an introduction to this research area). 

In this paper, we explore to what extent this attractive feature of sparse polynomial approximation methods can be leveraged in the setting of high-dimensional PDEs. As a model problem, we consider the high-dimensional diffusion equation. This serves as a simplified model for problems like the Black-Scholes, Schr\"odinger and Fokker-Planck equation, whose higher-order term is a second-order diffusion operator. Moreover, in order to keep the theoretical analysis accessible, we consider periodic boundary conditions or, equivalently, we assume the PDE domain to be the $d$-dimensional torus. 

\subsection{Main contributions}
\label{s:contributions}

We start by summarizing our main contributions.

\paragraph{Proposed method.} We propose a new numerical  method for the numerical solution of high-dimensional PDEs called compressive Fourier collocation, based on Monte Carlo sampling and sparse recovery via greedy algorithms or $\ell^1$ minimization. Our work builds upon the compressive spectral collocation method proposed in \cite{brugiapaglia2020compressivespectral}. We note, however, that the method in \cite{brugiapaglia2020compressivespectral} is heavily affected by the curse of dimensionality due to the use of (subsampled) tensorized equispaced grids and truncation sets of tensor product type, which lead to an exponential dependence of the number of collocation points on the dimension. We overcome this limitation thanks to the use of Monte Carlo sampling and truncation sets of hyperbolic cross type. %Notice that there is a key difference with respect to the compressive spectral collocation method proposed in \cite{brugiapaglia2020compressivespectral}: the method proposed in this paper entirely avoids the use of structured (subsampled) grids and only resorts to Monte Carlo sampling. This feature leads to lessening the curse of dimensionality in the number of collocation points.

\paragraph{Theoretical analysis.} We carry out a rigorous theoretical analysis of the proposed method for the solution of the diffusion equation on the $d$-dimensional torus. Our analysis shows that, under suitable regularity conditions on the diffusion coefficient, the compressive Fourier collocation method is \emph{provably able to lessen the curse of dimensionality with respect to the number of collocation points}. Specifically, it requires a number of collocation points that scales only logarithmically with the dimension $d$ of the PDE domain. This is proved in Theorem~\ref{thm:recovery_CFC}. 

Our analysis is based on the framework of random sampling in bounded Riesz system, recently introduced in \cite{brugiapaglia2021sparse}. The main technical core of our analysis is aimed at finding sufficient conditions on the diffusion coefficient that allow us to recast compressive Fourier collocation in this framework. We note that in our method the spectral basis  (i.e., a renormalized Fourier system) is orthogonal. However, our analysis relies on showing that the family of functions obtained by applying the PDE operator to the spectral basis is a bounded Riesz system under sufficient conditions on the PDE coefficients. This is proved in Propositions~\ref{prop:1-sparse_diffusion} and \ref{prop:compressible_eta}. 

\paragraph{Numerical experiments.} We numerically validate the proposed method and show that it is able to approximate both sparse and nonsparse solutions accurately. We show that the performance of the method is robust to variations in the diffusion coefficient and with respect to different parameter settings. Our numerical results involve simulations up to dimension $d=20$.

\subsection{Related literature} 
\label{s:literature}

The literature on numerical methods for high-dimensional PDEs is vast and fast-growing. Although a detailed review is outside the scope of this work, we provide a (noncomprehensive) set of references to some key papers in the field. 

A first class of methods for solving PDEs in moderately-high dimension is based on sparse grids; see, e.g., \cite{shen2010efficient,shen2012efficient}. A more recent line of work focused on constructing high-dimensional PDE solvers that are able to exploit low-rank structures in tensor-based approximations; see, e.g., \cite{bachmayr2015adaptive,bachmayr2016tensor,dahmen2016tensor}. In the past few years, a research direction that has gained considerable attention in the scientific machine learning community is the development of high-dimensional PDE solvers based on deep neural networks; see, e.g.,  \cite{berner2020analysis,elbrachter2022dnn,grohs2020deep,gu2021selectnet,han2018solving}. We refer to  \cite{weinan2021algorithms} for an extensive literature review on the topic.

From a different perspective, our work can also be contextualized within the literature of numerical methods for PDEs based on compressive sensing. Compressive (or compressed) sensing was introduced in \cite{donoho2006compressed,candes2006robust} and had a transformative impact on an interdisciplinary array of fields including, notably, scientific computing. For general introductions to compressive sensing and its applications, we refer to \cite{adcock2021compressive,foucart2013mathematical,lai2021sparse}. As mentioned in \S\ref{s:contributions}, our work builds upon the compressive spectral collocation method proposed in \cite{brugiapaglia2020compressivespectral}. Yet, we recall that the method proposed in our paper is provably able to lessen the curse of dimensionality with respect to the number of collocation points, as opposed to the method in \cite{brugiapaglia2020compressivespectral}, thanks to the improved truncation and collocation strategies.
Other numerical methods for PDEs based on compressive sensing include techniques based on Petrov-Galerkin discretizations  \citep{jokar2010sparse,brugiapaglia2015compressed, brugiapaglia2018theoretical, brugiapaglia2021wavelet} and  isogeometric analysis  \citep{brugiapaglia2020compressive, kang2019economical}. Another related line of work is the development of sparsity-promoting spectral methods for
multiscale problems based on sparse Fourier transforms \citep{daubechies2007sparse} and soft thresholding \citep{mackey2014compressive,schaeffer2013sparse}. In this direction, a very promising class of sublinear-time algorithms based on the sparse Fourier transform and on Fourier-Galerkin discretizations for high-dimensional, multiscale elliptic PDEs with periodic boundary conditions has been recently proposed in \cite{gross2023sparse}. For further pointers to the literature and historical remarks, we refer to \cite[Section 1.2]{brugiapaglia2020compressivespectral}.

Finally, we observe that using truncation sets of hyperbolic cross type in combination with Monte Carlo sampling is a well-known strategy in the context of high-dimensional \emph{polynomial} approximation and its applications to \emph{parametric} PDEs (see \cite{adcock2022sparse} and references therein). However, to the best of our knowledge, this is the first work where this strategy is successfully employed to construct a sparse high-dimensional PDE solver.

\subsection{Outline of the paper} 

The paper is organized as follows. \S\ref{s:setup} defines notation used throughout the paper, the model problem, and introduces the compressive Fourier collocation method. In \S\ref{s:theory}, we illustrate theoretical recovery guarantees for the compressive Fourier collocation method for the solution of the periodic diffusion equation. In \S\ref{s:numerics}, we present numerical experiments in dimension two and eight with different diffusion coefficients and for sparse and nonsparse exact solutions. In \S\ref{s:conclusions}, we conclude by summarizing our main findings and discussing avenues of future research. The proofs of the main results and further details on the numerical experiments are presented in Appendix~\ref{s:appendix}.

\section{Problem setting}
\label{s:setup}

In this section, we summarize standard mathematical notation used throughout the paper in \S\ref{s:notation}, define the high-dimensional periodic diffusion equation in \S\ref{s:model_problem}, and illustrate the compressive Fourier collocation method in \S\ref{s:CFC_method}. 

\subsection{Notation}\label{s:notation}

We denote the set of positive integers as $\mathbb{N}$, the set of nonnegative integers as $\mathbb{N}_0$ and the set of integers as $\mathbb{Z}$. For any $n \in \mathbb{N}$, the set of first $n$ nonnegative integers is denoted by $[n] = \{1,\ldots,n\}$. The notation $x \lesssim y$ means that there exists a universal constant $C >0$ such that $x \leq C y$. Moreover, when $X$ is a set, $|X|$ denotes its cardinality. The Kronecker delta is denoted by $\delta_{i,j}$. 

We equip $\mathbb{C}^n$ with the standard inner product $\bm{z}\cdot\bm{w} = \sum_{i\in[n]} z_i \bar{w}_i$ and its induced norm $\|\bm{z}\|_2 = (\bm{z} \cdot \bm{z})^{1/2}$. In general, we denote the $\ell^p$-norm of a vector $\bm{z} \in \mathbb{C}^n$ as $\|\bm{z}\|_p = (\sum_{i \in [n]} |z_i|^p)^{1/p}$, for any $1\leq p\leq \infty$, and its $\ell^0$-``norm'' as $\|\bm{z}\|_0 = |\supp(\bm{z})|$, where $\supp(\bm{z})=|\{j\in [N] : z_j \neq 0\}|$. These definitions naturally extend to the case where $\bm{z}$ is an infinite complex-valued sequence. Moreover, given a vector (or sequence)  $\bm{z} = (z_j)_{j \in \Lambda}$ indexed over a countable set $\Lambda$,  we define the restriction of $\bm{z}$ to $S$ as $\bm{z}_S = (z_j)_{j\in S}$ for any subset $S\subseteq \Lambda$. We denote the space of sequences $\bm{z} = (z_i)_{i\in \Lambda}$ such that $\|\bm{z}\|_p < \infty$ as $\ell^p(\Lambda;\mathbb{C})$. When $|\Lambda|=N <\infty$, we will identify $\ell^p(\Lambda;\mathbb{C}) \cong \mathbb{C}^N$ through the natural isomorphism.

We denote the one-dimensional torus by $\mathbb{T} = \mathbb{R}/\sim$, where $\sim$ is the equivalence relation on $\mathbb{R}$ defined by $x\sim y$ if and only if $x - y \in \mathbb{Z}$. We denote Lebesgue and Sobolev spaces over the $d$-dimensional torus as $L^p(\mathbb{T}^d)$ and $H^k(\mathbb{T}^d)$, with $1\leq p \leq \infty$ and $k \in \mathbb{N}_0$, with the convention that $H^0 = L^2$, and where we assume the functions to be complex-valued. The Lebsegue space $L^2$ is equipped with the inner product $\langle u, v \rangle = \int_{\mathbb{T}^d} u(\bm{x}) \overline{v(\bm{x})} \mathrm{d}\bm{x}$. We recall that the $H^1$-norm is given by $\|u\|_{H^1} = \|u\|_{L^2} + |u|_{H^1}$, where $\|u\|_{L^p}^p = \int_{\mathbb{T}^d}|u(\bm{x})|^p \mathrm{d}\bm{x}$ defines the $L^p$-norm for any $p\geq 1$ and $|u|_{H^1}^2 = \int_{\mathbb{T}^d}\|\nabla u(\bm{x})\|_2^2 \mathrm{d}\bm{x}$ is the $H^1$-seminorm.

\subsection{Model problem: The periodic diffusion equation}
\label{s:model_problem}

Our model problem is a diffusion equation over the torus $\mathbb{T}^d$, with $d \in \mathbb{N}$. Our interest in this paper lies in the scenario where $d \gg 1$. In addition to periodic boundary conditions, we add a zero-mean linear constraint to the equation in order for the problem to be well posed.  This leads to the following equations:
\begin{align}
& - \nabla \cdot (a(\bm{x}) \nabla u(\bm{x})) = f(\bm{x}), \quad \forall \bm{x} \in \mathbb{T}^d, \label{eq:diffusion_eq_periodic_1}\\
& \int_{\mathbb{T}^d} u(\bm{x}) \mathrm{d} \bm{x}= 0,\label{eq:diffusion_eq_periodic_2}
\end{align}
where the diffusion coefficient $a$ is such that 
\begin{equation}
\label{eq:suff_cond_diffusion}
a \in C^1(\mathbb{T}^d)
\quad \text{and} \quad
\min_{\bm{x} \in \mathbb{T}^d}\Real(a(\bm{x})) \geq a_{\min} >0,
\end{equation} 
in order to guarantee ellipticity, and the forcing term $f \in L^2(\mathbb{T}^d)$. Using the regularity theory of second-order elliptic problems (see, e.g., \cite[Section 6.3]{evans2010partial}), these assumptions on $a$ and $f$ guarantee that weak solutions to \eqref{eq:diffusion_eq_periodic_1}--\eqref{eq:diffusion_eq_periodic_2} satisfy $u\in H^2(\mathbb{T}^d)$ and lead to considering the following solution space: 
\begin{equation}
\label{eq:solution_space}
U = \left\{ v \in H^2 (\mathbb{T}^d): \int_{\mathbb{T}^d} v(\bm{x}) \mathrm{d} \bm{x} = 0 \right\}.
\end{equation}

\subsection{Compressive Fourier collocation} 
\label{s:CFC_method}

We now introduce the compressive Fourier collocation method.

\paragraph{Discretization of the PDE.}
Let us assume to have a basis $\{\Psi_{\bm{\nu}}\}_{\bm{\nu} \in \mathbb{Z}^d}$ (called the \emph{spectral} basis) for the solution space $U$ in \eqref{eq:solution_space}. This choice ensures that the linear constraint \eqref{eq:diffusion_eq_periodic_2} is enforced in a natural way. Then we can expand our solution as 
$$
u = \sum_{\bm{\nu} \in \mathbb{Z}^d} c_{\bm{\nu}} \Psi_{\bm{\nu}}.
$$
Since this expansion has infinitely many terms, it does not lead to an implementable approximation method. Hence, we consider a finite multi-index set $\Lambda \subset \mathbb{Z}^d$ with $|\Lambda| < \infty$. This leads to the finite expansion
$$
u_{\Lambda} = \sum_{\bm{\nu} \in \Lambda} c_{\bm{\nu}} \Psi_{\bm{\nu}}.
$$
We want to collocate the diffusion equation \eqref{eq:diffusion_eq_periodic_1} using the finite spectral basis $\{\Psi_{\bm{\nu}}\}_{\bm{\nu} \in \Lambda}$ by means of Monte Carlo sampling. Hence, we randomly generate $m$ independent points
$$
\bm{y}_1,\ldots,\bm{y}_m \in \mathbb{T}^d,
$$
uniformly distributed over $\mathbb{T}^d$. Letting $N=|\Lambda|$, we assume to have an ordering for the multi-indices in 
$\Lambda = \{\bm{\nu}_1, \ldots,\bm{\nu}_N\}$ (e.g., the lexicographic ordering).  This leads to the linear system
\begin{equation}
\label{eq:SC_system}
A \bm{z} = \bm{b},
\end{equation}
where $A \in \mathbb{C}^{m \times N}$ (called compressive Fourier collocation matrix) and $\bm{b} \in \mathbb{C}^m$ are defined by
\begin{equation}
\label{eq:def_A_b}
A_{ij}  = \frac{1}{\sqrt{m}}[-\nabla \cdot (a \nabla \Psi_{\bm{\nu}_j})](\bm{y}_i) 
\quad \text{and} \quad 
b_{i}  = \frac{1}{\sqrt{m}} f(\bm{y}_i), \quad \forall i \in [m], j \in [N].
\end{equation}
The normalization factor $1/\sqrt{m}$ is needed for technical reasons explained in Appendix~\ref{app:Gram}.

\paragraph{The spectral basis.}
Due to the presence of periodic boundary conditions, a natural choice for the spectral basis is the (complex) Fourier basis. The $L^2$-orthonormal Fourier basis is defined as
\begin{equation}
\label{eq:Fourier_system}
F_{\bm{\nu}}(\bm{x}) = \exp(2 \pi \mathrm{i} \, \bm{\nu} \cdot \bm{x}), \quad \forall \bm{\nu} \in \mathbb{Z}^d, \; \forall \bm{x} \in \mathbb{T}^d.
\end{equation}
The Fourier system $\{F_{\bm{\nu}}\}_{\bm{\nu} \in \mathbb{Z}^d}$ is a bounded orthonormal system of $L^2(\mathbb{T}^d)$. In fact, $\langle F_{\bm{\nu}}, F_{\bm{\mu}}\rangle = \delta_{\bm{\nu},\bm{\mu}}$ and $\|F_{\bm{\nu}}\|_{L^\infty} = 1$, for all $\bm{\nu},\bm{\mu} \in \mathbb{Z}^d$. However, since our spectral collocation matrix involves evaluations of the PDE operator applied to the basis functions, we need to consider a rescaling of this system. Let us consider, for the sake of simplicity, the case of the Poisson equation, where $a \equiv 1$. Applying the PDE operator to the $L^2$-normalized Fourier basis functions yields 
$-\Delta F_{\bm{\nu}} = 4\pi^2 \|\bm{\nu}\|_2^2 F_{\bm{\nu}}$, for all $\bm{\nu} \in \mathbb{Z}^d$. This leads us to define the spectral basis as the following rescaled version of the Fourier basis:
\begin{equation}
\label{eq:def_Psi}
\Psi_{\bm{\nu}} = \frac{1}{4\pi^2 \|\bm{\nu}\|_2^2} F_{\bm{\nu}}, \quad \forall \bm{\nu} \in \mathbb{Z}^d.
\end{equation}
This implies that $\{-\Delta \Psi_{\bm{\nu}}\}_{\bm{\nu} \in \mathbb{Z}^d} = \{F_{\bm{\nu}}\}_{\bm{\nu} \in \mathbb{Z}^d}$ is a bounded orthonormal system of $L^2(\mathbb{T}^d)$. At least in the case $a \equiv 1$, this puts us in the optimal position to perform compressive collocation, as the spectral collocation matrix $A$ is the sampling matrix associated with the random sampling in a bounded orthonormal system and it is therefore an ideal scenario for compressive sensing (see, e.g., \citet{foucart2013mathematical}). 
A substantial portion of our efforts in this paper will be devoted to showing that this normalization choice is appropriate for performing compressive sensing also in the case $a \not\equiv 1$. In essence, this is due to the fact that, under suitable sufficient conditions on the diffusion coefficient $a$, the system $\{-\nabla\cdot (a \nabla \Psi_{\bm{\nu}})\}_{\bm{\nu} \in \mathbb{Z}^d}$ is a \emph{bounded Riesz system} in $L^2(\mathbb{T}^d)$, following the framework of \cite{brugiapaglia2021sparse}. This will be discussed in detail in \S\ref{s:theory}. %\purple{[Point out that this is a form of diagonal preconditioning.]}

\begin{figure}[!t]%
\begin{subfigure}{.5\textwidth}
    \centering
    \includegraphics[width=\linewidth]{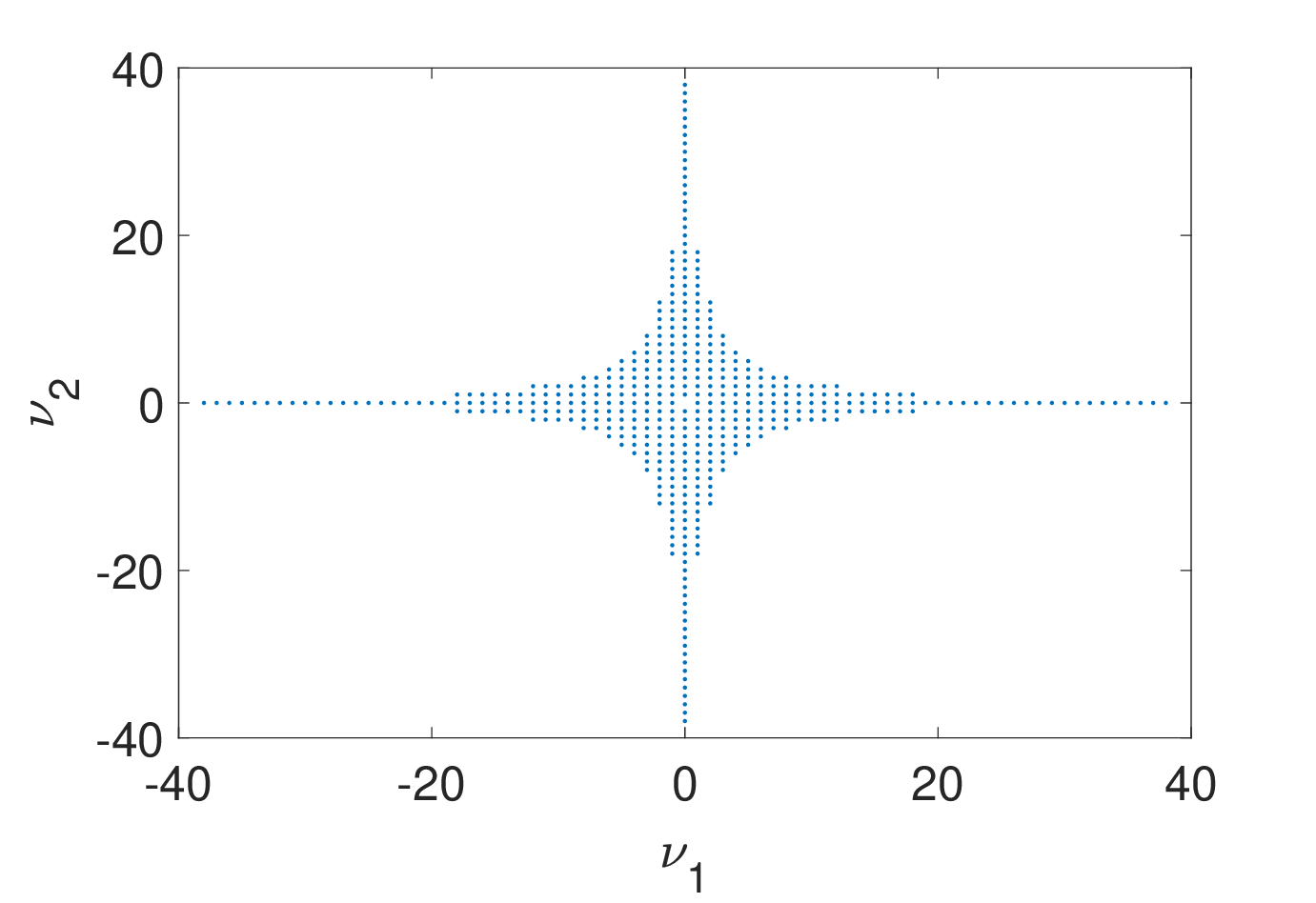}
    \caption{$d = 2$, $n=39$}
\end{subfigure}
\begin{subfigure}{.5\textwidth}
    \centering
    \includegraphics[width=\linewidth]{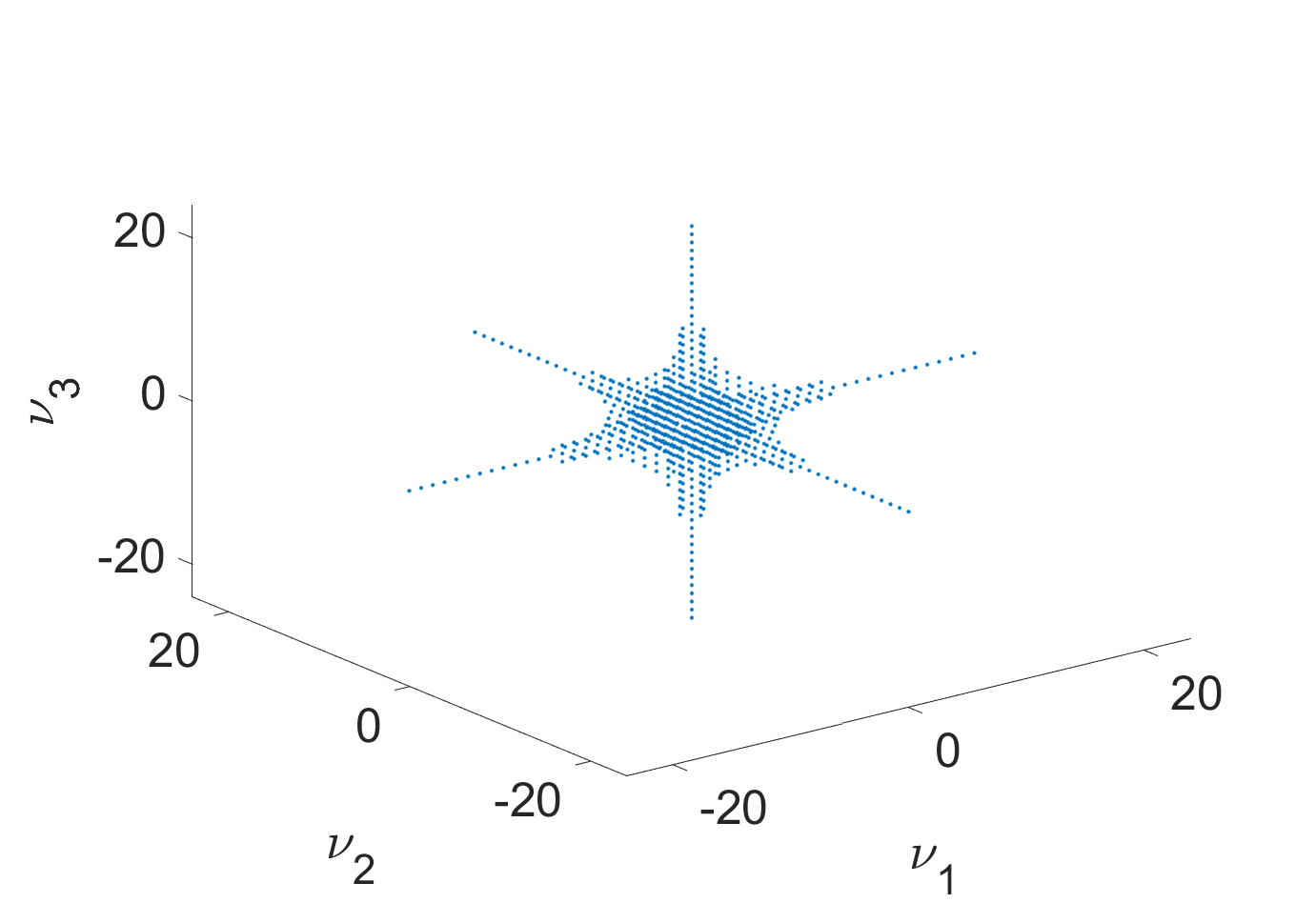}
    \caption{$d = 3$, $n=25$}
\end{subfigure}
\caption{\label{fig:HC}The hyperbolic cross index set (minus the zero multi-index) $\Lambda$ defined in \eqref{eq:def_Lambda} for different values of $d$ and $n$.}
\end{figure}

As for the truncation set $\Lambda \subset \mathbb{Z}^d$, we will choose a hyperbolic cross of $\mathbb{Z}^d$ (minus the zero multi-index), i.e.\
\begin{equation}
\label{eq:def_Lambda}
\Lambda = \left\{\bm{\nu}\in \mathbb{Z}^d : \prod_{l = 1}^d (|\nu_l| +1) \leq n  \right\} \setminus \{\bm{0}\}
\end{equation}
Fig.~\ref{fig:HC} shows $\Lambda$ in dimension $d=2$ and $d=3$. 
The main advantage of this choice is that the cardinality of the hyperbolic cross has a controlled growth with respect to $d$, as opposed to other types of multi-index set such as tensor product or total degree index sets (see, e.g., \cite{adcock2022sparse}). Specifically, combining cardinality bounds proved in \cite{chernov2016new,kuhn2015approximation}, it is possible to see that
\begin{equation}
\label{eq:cardinality_bound_Lambda}
N = |\Lambda| \leq \min\{4 n^5 16^d, \mathrm{e}^2 n^{2+\log_2(d)}\}.
\end{equation}
For more details, see Appendix~\ref{app:HCbounds}. This cardinality bound will be crucial to prove that compressive Fourier collocation is able to lessen the curse of dimensionality with respect to the number of collocation points. Moreover, hyperbolic crosses are also well suited for effective high-dimensional approximation. In fact, Fourier approximations supported on hyperbolic crosses are known to be accurate when the target function is periodic and has mixed smoothness. For further details, we refer to \cite{dung2018hyperbolic,temlyakov2018multivariate}. 

\paragraph{Solving the linear system.}
In order to compute an approximation to the solution coefficients $\bm{c}_\Lambda=(c_{\bm{\nu}})_{\bm{\nu} \in \Lambda}$, we need to compute an (approximate) solution $\hat{\bm{c}}=(\hat{c}_{\bm{\nu}})_{\bm{\nu} \in \Lambda} \in \mathbb{C}^N$ to the linear system \eqref{eq:SC_system}. %To do this, we will assume the existence of a \emph{decoder} (borrowing a terminology that is very common in compressed sensing \cite{cohen2009compressed})
%$$
%\mathcal{D}: \mathbb{C}^{m\times N} \times \mathbb{C}^m \times \Theta \to \mathbb{C}^N.
%$$
%Given a matrix $A \in \mathbb{C}^{m \times N}$, a vector $\bm{b} \in \mathbb{C}^m$ and a set of tuning parameters $\theta \in \Theta$, the decoder $\mathcal{D}$ produces an approximation $\hat{\bm{c}}$ to the linear system \eqref{eq:SC_system}, namely
%$$
%\hat{\bm{c}} = \mathcal{D}(A, \bm{b}, \theta) \in \mathbb{C}^N,
%$$
%such that $A \hat{\bm{c}} \approx \bm{b}$. In practice, $\mathcal{D}$ will either rely on the numerical solution of a certain (convex) optimization problem or will be the output of a certain algorithm. 
The approximation to $u$ associated with the coefficients $\hat{\bm{c}}$ is defined by
$$
\hat{u} = \sum_{\bm{\nu} \in \Lambda} \hat{c}_{\bm{\nu}} \Psi_{\bm{\nu}}.
$$
%We have several options, depending on the fact that the system under exam is \emph{overdetermined} ($m \geq N$) or \emph{underdetermined} ($m < N$).
When $m \geq N$, the most natural way to find an approximate solution to \eqref{eq:SC_system} is via ordinary least squares by letting
\begin{equation}
\label{eq:OLS}
\hat{\bm{c}} \in \arg\min_{\bm{z} \in \mathbb{C}^N} \|A \bm{z} - \bm{b}\|_2^2.
\end{equation}
(Note that in cases where the minimizer is not unique we assume $\hat{\bm{c}}$ to be \emph{any} solution to the optimization problem of interest). However, due to the fact that $d\gg 1$, we would like to use less than $N$ collocation points in order to avoid the curse of dimensionality.

When $m < N$, the system \eqref{eq:SC_system} admits infinitely many solutions and we need to introduce some form of regularization or, equivalently, a structural \emph{a priori} assumption on the solution $u$ that allows us to retrieve an accurate approximation of it from the underdetermined linear system. In this paper, we will achieve this by assuming that the solution $u$ be \emph{sparse} or \emph{compressible} with respect to the spectral basis $\{\Psi_{\bm{\nu}}\}_{\bm{\nu} \in \mathbb{Z}^d}$. We say that a solution $u = \sum_{\bm{\nu} \in\mathbb{Z}^d} c_{\bm{\nu}} \Psi_{\bm{\nu}}$ is $s$-sparse (with respect to the system $\{\Psi_{\bm{\nu}}\}_{\bm{\nu} \in \mathbb{Z}^d}$) if $\|\bm{c}\|_0 \leq s$. Moreover, informally, $u$ is compressible if its best $s$-term approximation error
$$
\sigma_s(\bm{c})_p = \inf_{\bm{z} \in \ell^p(\mathbb{Z}^d)} \|\bm{c} - \bm{z}\|_p,
$$
has a fast decay rate with respect to $s$ for some $p\geq 1$. Rigorous decay rates for the best $s$-term approximation error of functions with mixed smoothness can be found in, e.g., \citet[Chapter 9]{temlyakov2018multivariate}. For example, \cite[Theorem~9.1.4]{temlyakov2018multivariate} implies bounds of the form $\sigma_s(\bm{c})_2 \lesssim m^{-r+\alpha}(\log m)^{(d-1)(r-\beta)}$ when $u$ is such that its partial derivatives $u^{(r_1,\ldots,r_d)}$ satisfy
$\|u^{(r_1,\ldots,r_d)}\|_{L^q} \leq 1$
for every $r_j \leq r \in \mathbb{N}$, where $\alpha,\beta \in [0, r)$ are suitable constants depending on $q > 1$. We also note that, under assumption \eqref{eq:suff_cond_diffusion} on the diffusion coefficient $a$, such mixed regularity conditions for the solution $u$ are implied by the analogous conditions on the forcing term $f$ (see, e.g., \cite[Section~6]{bungartz1999note} and  \cite{griebel2009optimized}). In this paper, we do not focus on any function classes in particular, but only assume that the solution $u$ is sparse or compressible. %Combining the recovery error estimates results of this paper with rigorous approximation theoretical results such as those in \cite{temlyakov2018multivariate} is out of the scope of this paper and it left as an important direction of future investigation (see also \S\ref{s:conclusions}).  

This leads us to employ sparse recovery techniques to compute approximate solutions to the underdetermined linear system \eqref{eq:SC_system}. We consider two approaches for this. The first one is recovery via Orthogonal Matching Pursuit (OMP). One of the attractive features of OMP is its computational efficiency for small values of the target sparsity.

\begin{algorithm}
\caption{Orthogonal Matching Pursuit (OMP)} \label{alg:OMP}
\hspace*{\algorithmicindent} \textbf{Input:} Measurement matrix $A\in\mathbb{C}^{m \times N}$, vector $\bm{b} \in \mathbb{C}^m$, number of iterations $K \in \mathbb{N}$. \\
\hspace*{\algorithmicindent} \textbf{Output:} A $K$-sparse vector $\hat{\bm{c}} \in \mathbb{C}^N $. 
\begin{algorithmic}[1]
\State Normalize columns of $A$, i.e. $A_{ij} \gets A_{ij}/\sqrt{\sum^m_{i=1}|A_{ij}|^2}$, $\forall i \in [m], \forall j \in [N]$ 
\State $S_0 \gets \emptyset,\ \bm{z}_0=\bm{0},\ n=0$\; 
\For {$n=0,\dots,K-1$}
    \State $j_{n+1} \gets \argmax_{j \in [N]} \left\{ \left|(A^*(\bm{b}-A \bm{z}_n))_j\right| \right\}$ 
    \State $S_{n+1} \gets S_{n} \cup j_{n+1}$
    \State $\bm{z}_{n+1} \gets \argmin_{\bm{z}\in \mathbb{C}^N} \left\{ \left\| \bm{b} - A \bm{z}\right\|_2,\,\mbox{supp}(\bm{z}) \subseteq S_{n+1} \right\}$
\EndFor
\State $\hat{\bm{c}} \gets \bm{z}_K$ 
\end{algorithmic}
\end{algorithm}

The second approach is based on $\ell^1$ minimization, via the Quadratically-Constrained Basis Pursuit (QCBP) convex optimization program
\begin{equation}
\label{eq:QCBP}
\hat{\bm{c}} \in \arg\min_{\bm{z} \in \mathbb{C}^N} \|\bm{z}\|_1 \quad \text{such that} \quad \|A \bm{z} - \bm{b}\|_2 \leq \eta,
\end{equation}
where $\eta > 0$ is a tuning parameter. We note in passing that other choices, such as the square-root LASSO convex optimization program, might be considered to improve robustness with respect to the tuning parameter choice (see, e.g., \cite{adcock2019correcting}). For the sake of clarity, we summarize the compressive Fourier collocation method in Algorithm~\ref{alg:Compressive_Fourier_collocation}.

\begin{algorithm}
\caption{Compressive Fourier collocation method for periodic diffusion equations} \label{alg:Compressive_Fourier_collocation}
\hspace*{\algorithmicindent} \textbf{Input:} Dimension $d\in \mathbb{N}$, diffusion coefficient $a$, forcing term $f$, hyperbolic cross order $n\in \mathbb{N}$, number of collocation points $m\in \mathbb{N}$. \\ 
\hspace*{\algorithmicindent} \textbf{Output:} Approximate solution $\hat{u}$.
\begin{algorithmic}[1]
\State Generate hyperbolic cross set $\Lambda=\{\bm{\nu}_1,\ldots,\bm{\nu}_m\}$ as in \eqref{eq:def_Lambda} 
\State Randomly generate $m$ independent collocation points $\bm{y}_1,\ldots,\bm{y}_m$ uniformly distributed in $\mathbb{T}^d$
\State Using the generated collocation points, build $A \in \mathbb{C}^{m \times N}$ and $\bm{b} \in \mathbb{C}^m$ as in \eqref{eq:def_A_b}
\State Compute approximate coefficients $\hat{\bm{c}}\in\mathbb{C}^N$ via OMP (Algorithm \ref{alg:OMP})  or QCBP (see \eqref{eq:QCBP})\;
\State Define $\hat{u} \gets \sum_{\bm{\nu} \in \Lambda} \hat{c}_{\bm{\nu}} \Psi_{\bm{\nu}}$
\end{algorithmic}
\end{algorithm}

\paragraph{Truncation error.} Note that, in general, $\bm{c}_\Lambda$ does not solve \eqref{eq:SC_system} exactly, but only approximately. In fact, the finite expansion $u_{\Lambda}$ of $u$ satisfies the equation
\begin{align}
& - \nabla \cdot (a(\bm{x}) \nabla u_\Lambda(\bm{x})) = \tilde{f}_{\Lambda}(\bm{x}), \quad \forall \bm{x} \in \mathbb{T}^d, \label{eq:diffusion_eq_periodic_finite_1}\\
& \int_{\mathbb{T}^d} u_\Lambda(\bm{x}) \,\mathrm{d} \bm{x}= 0, \label{eq:diffusion_eq_periodic_finite_2}
\end{align}
where $\tilde{f}_{\Lambda}(\bm{x}) = f(\bm{x}) + \nabla \cdot (a(\bm{x}) \nabla (u(\bm{x}) - u_{\Lambda}(\bm{x}))$.
Therefore, 
\begin{equation}
\label{eq:linear_system_c_Lambda}
A \bm{c}_{\Lambda} = \bm{b} + \bm{e},
\end{equation}
where the error vector $\bm{e} \in \mathbb{C}^m$ is defined by
\begin{equation}
    \label{eq:def_error_vector}
    e_{i} = \frac{1}{\sqrt{m}}[\nabla \cdot (a \nabla (u - u_{\Lambda}))](\bm{y}_i),
    \quad \forall i \in [m].
\end{equation} 
Note that, in general,
\begin{equation}
\label{eq:error_norm_bound}
\|\bm{e}\|_2 \leq \|\nabla \cdot (a \nabla (u - u_{\Lambda}))\|_{L^\infty}.
\end{equation}
This term measures the truncation error (due to the introduction of a finite set $\Lambda \subset \mathbb{Z}^d$) in a metric that depends on the operator defining the PDE and, in particular, on the diffusion coefficient $a$. Note also that, thanks to the strong law of large numbers, we have
$$
\|\bm{e}\|_2^2
= \frac{1}{m} \sum_{i = 1}^m |\nabla \cdot (a\nabla (u-u_{\Lambda}))(\bm{y}_i)|^2
\to \|\nabla \cdot (a\nabla (u-u_{\Lambda}))\|_{L^2}^2,\quad \text{as $m \to \infty$,}
$$
with convergence rate $1/\sqrt{m}$. However, in the following, we will use the upper bound \eqref{eq:error_norm_bound} in order to avoid further technical difficulties due to the asymptotic nature of this estimate.

\section{Theoretical analysis}
\label{s:theory}
We now present rigorous theoretical recovery guarantees for the compressive Fourier collocation method. Our analysis relies on the framework of random sampling in bounded Riesz systems recently introduced in \cite{brugiapaglia2021sparse}. Our main technical effort is devoted to showing that the compressive Fourier collocation matrix $A$ defined in \eqref{eq:def_A_b} is obtained by sampling a bounded Riesz system. This is discussed in \S\ref{s:bounded_Riesz}. After doing this, in \S\ref{s:recovery} we apply sparse recovery results for random sampling in bounded Riesz systems to our case and prove error bounds for the compressive Fourier collocation solutions computed via OMP and QCBP. The proofs of the results stated in this section can be found in Appendix~\ref{s:appendix}.

\subsection{Compressive Fourier collocation and random sampling in bounded Riesz systems}
\label{s:bounded_Riesz}
We start by recalling the definition of the bounded Riesz system. We restrict our attention to Riesz systems in $L^2(\mathbb{T}^d)$, although the definition can be extended to general complex Hilbert spaces.

\begin{definition}[Bounded Riesz System]
\label{def:Riesz}
Let $\Lambda \subset \mathbb{Z}^d$ be a countable set and let $0<b_\Phi \leq B_\Phi < \infty$. A set of functions $\{\Phi_{\bm{\nu}}\}_{\bm{\nu} \in \Lambda} \subset L^2(\mathbb{T}^d)$ is a \emph{Riesz system} with constants $b_{\Phi}$ and $B_{\Phi}$ if 
$$
b_\Phi\|\bm{z}\|^2_2 
\leq \left\| \sum_{\bm{\nu} \in \Lambda} z_{\bm{\nu}}
\Phi_{\bm{\nu}}\right\|^2_{L^2} 
\leq B_\Phi \|\bm{z}\|^2_2, 
\quad \forall \bm{z}=(z_{\bm{\nu}})_{\bm{\nu} \in \Lambda} \in \ell^2(\Lambda;\mathbb{C}).
$$
The constants $b_{\Phi}$ and $B_{\Phi}$ are called \emph{lower and upper Riesz constants}, respectively. Moreover, the system $\{\Phi_{\bm{\nu}}\}_{\bm{\nu} \in \Lambda}$ is \emph{bounded} if there exists a constant $0< K_{\Phi} < \infty$ such that
$$
\|\Phi_{\bm{\nu}}\|_{L^\infty} =\esssup_{\bm{x} \in \mathbb{T}^d} |\Phi_{\bm{\nu}}(\bm{x})| \leq K_{\Phi}, \quad \forall \bm{\nu} \in \Lambda.
$$
\end{definition}
Note that any $L^2$-orthonormal system is a Riesz system with $b_{\Phi} = B_{\Phi} = 1$. In particular, the Fourier system $\{F_{\bm{\nu}}\}_{\bm{\nu}\in\mathbb{Z}^d}$ is a bounded Riesz system with $b_{\Phi} = B_{\Phi} = K_{\Phi} = 1$

We define the second-order operator $\mathscr{L}:H^2(\mathbb{T}^d) \to L^2(\mathbb{T}^d)$ associated with \eqref{eq:diffusion_eq_periodic_1} as 
$$
\mathscr{L}[v] = -\nabla \cdot (a \nabla v), \quad \forall v \in H^2(\mathbb{T}^d).
$$
Using this notation, we also define 
\begin{equation}
\label{eq:def_Phi}
\Phi_{\bm{\nu}} = \mathscr{L}[\Psi_{\bm{\nu}}], \quad \forall  \bm{\nu} \in \Lambda.
\end{equation}
Recall that we are assuming $\bm{0} \notin \Lambda$ in order to naturally satisfy the condition \eqref{eq:diffusion_eq_periodic_2}.
As a consequence of these definitions, the compressive Fourier collocation matrix $A$ defined in \eqref{eq:def_A_b} admits the equivalent definition $A_{ij} = \Phi_{\bm{\nu}_j}(\bm{x}_i)$, for all $i \in [m]$ and $j \in [N]$ and is hence a random sampling matrix associated with the system $\{\Phi_{\bm{\nu}}\}_{\bm{\nu} \in \Lambda}$. For this reason, showing that $\{\Phi_{\bm{\nu}}\}_{\bm{\nu} \in \Lambda}$ is a bounded Riesz system will frame the compressive Fourier collocation problem in the setting of \cite{brugiapaglia2021sparse} and allow us to apply the corresponding sparse recovery theorems to quantify the accuracy of approximations computed via OMP or QCBP. 

Note that the operator $\mathscr{L}$, and hence the system $\{\Phi_{\bm{\nu}}\}_{\bm{\nu} \in \Lambda}$ depend on the diffusion coefficient $a$. Hence, as explained in detail in Appendix~\ref{app:Gram}, we need to track the dependence of the Gram matrix $G$ of $\{\Phi_{\bm{\nu}}\}_{\bm{\nu}\in\Lambda}$ on $a$ in order to study the Riesz constants. Fortunately, it is possible to obtain an explicit formula for the entries of the Gram matrix $G$ in terms of the Fourier coefficients of the diffusion coefficient $a$ (see Lemma~\ref{lem:expansion_G} in Appendix~\ref{app:Gram}). For this reason, we let
\begin{equation}
\label{eq:eta_expansion}
a = \sum_{\bm{\nu} \in \mathbb{Z}^d} e_{\bm{\nu}} F_{\bm{\nu}}, \quad \text{with } e_{\bm{\nu}} = \langle a, F_{\bm{\nu}}\rangle, \quad \forall \bm{\nu} \in \mathbb{Z}^d,
\end{equation} 
and where $\{F_{\bm{\nu}}\}_{\bm{\nu} \in \mathbb{Z}^d}$ is the $L^2$-orthonormal Fourier system defined in \eqref{eq:Fourier_system}. Using the explicit formula for the Gram matrix $G$ in terms of the coefficients $\bm{e}$ of $a$ (provided by Lemma~\ref{lem:expansion_G} in Appendix~\ref{app:Gram}) and the Gershgorin circle theorem (recalled in Appendix~\ref{app:Gershgorin}), we can find explicit bounds for the minimum and maximum eigenvalues of $G$, which imply estimates for the lower and upper Riesz constant $b_{\Phi}$ and $B_{\Phi}$ of $\{\Phi_{\bm{\nu}}\}_{\bm{\nu}\in\Lambda}$, under suitable sufficient conditions on the diffusion coefficient $a$. We illustrate this in two examples of increasing generality (and technical complexity).

 We first analyze a very simple model case where the diffusion term is composed by a 1-sparse perturbation of a constant. %In other words, we assume $a$ to be composed by the superposition of a constant function and of a multiple of a Fourier basis function $F_{\bm{\nu}^*}$. 
 The following result is proved in Appendix~\ref{app:prop:1-sparse_diffusion}. %\purple{[In this simple example we could just set $e_{\bm{0}} = 1$ to have an easier statement and to claim that this is a `` 1-sparse perturbation of the Poisson problem''.]}
\begin{proposition}[Bounded Riesz property: $1$-sparse perturbation of constant diffusion]
\label{prop:1-sparse_diffusion}
Consider a diffusion coefficient of the form $a = e_{\bm{0}}+e_{\bm{\nu}^*}F_{\bm{\nu}^*}$, for some $\bm{\nu}^* \in \mathbb{Z}^d$ and with $e_{\bm{0}},e_{\bm{\nu}^*}\in \mathbb{C}$ such that $\alpha := |e_{\bm{\nu}^*}| (2\|\bm{\nu}^*\|_2+3)<|e_{\bm{0}}|$. Then, for any $\Lambda \subset \mathbb{Z}^d$, with $\bm{0} \notin \Lambda$, the system $\{\Phi_{\bm{\nu}}\}_{\bm{\nu} \in \Lambda}$ defined in \eqref{eq:def_Phi} is a bounded Riesz system in the sense of Definition~\ref{def:Riesz} with constants
$$
b_{\Phi}= \left|e_{\bm{0}}\right|^2- \left|e_{\bm{0}}\right| \alpha > 0, \quad  
B_{\Phi} = \left(|e_{\bm{0}}| + \frac{\alpha}{2}\right)^2, \quad \text{and} \quad 
K_{\Phi} = \left| e_{\bm{0}} \right| + \frac{\alpha}{2}.
$$

%\purple{[To do: Make this statement generic, with assumption of the form $|c_{\bm{\nu^*}}|<\frac{\xi}{\|\bm{\nu^*}\|_2+3}$.]}
\end{proposition}
Despite its simplicity, Proposition~\ref{prop:1-sparse_diffusion} illustrates an important point: diffusion coefficients that are highly oscillatory lead to worse Riesz constants. % and, hence, makes the problem more challenging (at least from a theoretical perspective). %In fact, a larger amplitude or frequency of the $1$-sparse perturbatio of $a$ leads to a smaller values of $b_{\Phi}$ and larger values of $B_{\Phi}$ and $K_{\Phi}$. 
Next, we consider a more general and challenging case. We study the case of a diffusion coefficient $a$ that is approximately sparse with respect to the Fourier basis. The following result highlights the impact of the approximate sparsity and the sparse approximation error of $a$ on the Riesz constants. This result is proved in Appendix~\ref{app:prop:compressible_eta}. 

\begin{proposition}[Bounded Riesz property: nonsparse diffusion coefficient]
\label{prop:compressible_eta}
Let $t\in \mathbb{N}$ and consider a diffusion coefficient $a$ satisfying \eqref{eq:suff_cond_diffusion} and with Fourier expansion \eqref{eq:eta_expansion}, with $e_{\bm{0}} \in \mathbb{R}$ and $e_{\bm{0}} \geq 0$, of the form 
\begin{equation*}
    a= a_t+a^*,
\end{equation*}
where $a_t$ is a $(t+1)$-sparse approximation of $a$ of the form
$$
a_t = e_{\bm{0}} + \sum_{\bm{\nu} \in T} e_{\bm{\nu}} F_{\bm{\nu}}, \quad \text{with } T \subset \mathbb{Z}^d \setminus \{\bm{0}\} \text{ and } |T|\leq t,
$$
and $a^* = a - a_t$ is the reminder term. Moreover, let $\Lambda \subset \mathbb{Z}^d \setminus \{\bm{0}\}$ with  $N = |\Lambda|$ and assume that
\begin{align}
    \beta & := \sqrt{t\left(\left\|a_t \right\|_{H^1}^2 - |e_{\bm{0}}|^2\right)} < (\sqrt{2}-1)|e_{\bm{0}}|, \label{eq:compressible_eta_cond_1}\\
    \gamma & :=\frac{\sqrt{N}}{2\pi} |a^*|_{H^1} 
    + \left\| a^* \right\|_{L^2} \leq  \sqrt{ |e_{\bm{0}}|^2- 2|e_{\bm{0}}| \beta - \beta^2 },\label{eq:compressible_eta_cond_2}
\end{align}
Then, the system $\{\Phi_{\bm{\nu}}\}_{\bm{\nu} \in \Lambda}$ defined in \eqref{eq:def_Phi} is a bounded Riesz system in the sense of Definition~\ref{def:Riesz} with constants
\begin{equation}
\label{eq:Riesz_constants_eta_compressible}
b_{\Phi} = \left( \sqrt{ |e_{\bm{0}}|^2- 2|e_{\bm{0}}| \beta - \beta^2}  -\gamma \right)^2>0, \quad  
B_{\Phi} = \left( \sqrt{ \left\|a_t \right\|_{H^1}^2 + 2|e_{\bm{0}}| \beta+ \beta^2} +\gamma \right)^2,
\end{equation}
and
\begin{equation}
\label{eq:unif_bound_eta_compressible}
K_{\Phi}  = |e_{\bm{0}}| + \beta + \left\| a^* \right\|_{L^\infty} + \sum_{l=1}^d \left\| \frac{\partial a^*}{\partial x_l}  \right\|_{L^\infty}.
\end{equation}
\end{proposition}
We observe that the two sufficient conditions \eqref{eq:compressible_eta_cond_1}--\eqref{eq:compressible_eta_cond_2} of Proposition~\ref{prop:compressible_eta} play substantially different roles. On the one hand, condition \eqref{eq:compressible_eta_cond_2} is a condition on the compressibility of the diffusion coefficient $a$ since it only involves the tail $a^*$ (and the constant term $e_{\bm{0}}$). On the other hand, condition \eqref{eq:compressible_eta_cond_1} only involves the $(t+1)$-sparse approximation $a_t$ of $a$. We observe that the constant $\gamma$ in \eqref{eq:compressible_eta_cond_2} depends on the truncation set $\Lambda$ through the term $\sqrt{N} = \sqrt{|\Lambda|}$. We suspect this might be an artifact of our proof strategy and understanding whether and how this term can be removed is an open problem. To provide further insights on conditions \eqref{eq:compressible_eta_cond_1}--\eqref{eq:compressible_eta_cond_2}, we illustrate a simple class of diffusion coefficients satisfying these assumptions in the following remark.
\begin{remark}[Example of diffusion coefficients satisfying the assumptions of Proposition~\ref{prop:compressible_eta}]
\label{app:example_diffusion}
Consider a non-negative multi-index $\bm{k} \in \mathbb{N}_0^d \setminus \{\bm{0}\}$, let $c_{\bm{0}},c_{\bm{k}} \in \mathbb{R}$, and define a real-valued diffusion coefficient of the form
\begin{equation}
\label{eq:exa_diffusion_cosine}
a(\bm{x}) = c_{\bm{0}} + c_{\bm{k}}\prod_{l=1}^d \cos{(2 \pi k_l x_l)}, \quad \forall \bm{x} \in \mathbb{T}^d,
\end{equation}
where $c_{\bm{0}}$ and $c_{\bm{k}}$ are such that $c_{\bm{0}} > |c_{\bm{k}}|$ (in order for condition \eqref{eq:suff_cond_diffusion} to be satisfied). This diffusion coefficient is 2-sparse with respect to the tensorized cosine basis $\{\prod_{l =1}^d\cos(2\pi k_l x_l)\}_{\bm{k} \in \mathbb{N}_0^d}$. However, it is $2^{\|\bm{k}\|_0}$-sparse with respect to the complex Fourier basis $\{F_{\bm{\nu}}\}_{\bm{\nu} \in \mathbb{Z}^d}$. In fact, for any $\bm{x} \in \mathbb{T}^d$, we have
\begin{align*}
    a(\bm{x}) & = c_{\bm{0}} + c_{\bm{k}}\prod_{l \in \supp(\bm{k})} \cos{(2 k_l \pi x_l)} \\
    & = c_{\bm{0}} + c_{\bm{k}}\prod_{l \in \supp(\bm{k}) } \frac{1}{2} \left( \exp{(2 k_l \pi \mathrm{i} x_l)} + \exp{(-2 k_l \pi \mathrm{i} x_l)} \right) \\
    & = c_{\bm{0}} + \frac{c_{\bm{k}}}{2^{\|\bm{k}\|_0}} \sum_{\substack{\bm{j} \in \{0,1\}^d\\ \supp(\bm{j}) \subseteq \supp(\bm{k})}}  F_{\bm{k} \odot (-1)^{\bm{j}}}(\bm{x}),
\end{align*}
where $\odot$ is the Hadamard (or componentwise) product operation and where we have denoted $(-1)^{\bm{j}}=((-1)^{j_1},(-1)^{j_2},\dots,(-1)^{j_d})$, for any $\bm{j} 
\in \{0,1\}^d$. 

Now, consider Proposition~\ref{prop:compressible_eta} with $t=2^{\|\bm{k}\|_0}$ and $T = \{\bm{k}\odot (-1)^{\bm{j}} : \bm{j} \in \{0,1\}^d \text{ and }\supp(\bm{j}) \subseteq \supp(\bm{k})\}$. First, condition \eqref{eq:compressible_eta_cond_2} is trivially satisfied since $a^* = 0$. Moreover, we see that 
\begin{align*}
    \beta & = \sqrt{t\left(\left\|a_t \right\|_{H^1}^2 - |e_{\bm{0}}|^2\right)} \\
    & = \sqrt{2^{\|\bm{k}\|_0} \sum_{\bm{\tau}\in T }\left( 1 + (2\pi)^2 \|\bm{\tau}\|_2^2\right)|e_{\bm{\tau}}|^2} \\
    & = \sqrt{2^{2\|\bm{k}\|_0} \left( 1 + (2\pi)^2 \|\bm{k}\|_2^2\right)\left(\frac{|c_{\bm{k}}|}{2^{\|\bm{k}\|_0}}\right)^2} \\
    & = |c_{\bm{k}}| \sqrt{1 + (2\pi)^2 \|\bm{k}\|_2^2}.
\end{align*}
Therefore, if $|c_{\bm{0}}| >  \frac{ \sqrt{1 + (2\pi)^2 \|\bm{k}\|_2^2}}{\sqrt{2}-1} |c_{\bm{k}}|$, then the diffusion coefficient $a$ in \eqref{eq:exa_diffusion_cosine} satisfies condition \eqref{eq:compressible_eta_cond_1} in Proposition \ref{prop:compressible_eta}. 
\end{remark}
Now that we have sufficient conditions on the diffusion coefficient $a$ in order for $\{\Phi_{\bm{\nu}}\}_{\bm{\nu} \in \Lambda}$ to be a bounded Riesz system, we are in a position to illustrate recovery guarantees for compressive Fourier collocation via QCBP and OMP. 

\subsection{Recovery guarantees compressive Fourier collocation}
\label{s:recovery}

%UNder the RIP, we have accurate and stable recovery guarantees for OMP and $\ell^1$ minimization (see recovery theorems in the Bounded Riesz basis paper). In the context of CSC, this implies... (add what was described in the previous section in terms of recovery error). So, Theorem [ref] implies the following.

%For QCBP we need consider \cite[Theorem 2.1]{cai2013sparse} with $t = 2$. This leads to $RIP(2s, \delta)$ with $\delta < 1/\sqrt{2})$ implies accurate and stable recovery via QCBP with usual bound and constants depending on delta.

Accurate and stable recovery guarantees for the compressive Fourier collocation method via OMP and QCBP are based on the framework of random sampling in bounded Riesz systems of \cite{brugiapaglia2021sparse}. The main elements of this framework are presented in Appendix~\ref{app:RIP}. Combining these with Proposition~\ref{prop:compressible_eta}  leads to the following recovery guarantees under sufficient regularity conditions on the diffusion coefficient $a$. The proof of this result can be found in Appendix~\ref{app:recovery_CFC}.

\begin{theorem}[Accurate and stable recovery for  compressive Fourier collocation]
\label{thm:recovery_CFC}
There exist universal constants $c, C_1, C_2>0$ such that approximations to the high-dimensional periodic diffusion equation \eqref{eq:diffusion_eq_periodic_1}--\eqref{eq:diffusion_eq_periodic_2} obtained via compressive Fourier collocation and computed via OMP (Algorithm~\ref{alg:OMP}) or QCBP (program \eqref{eq:QCBP}) satisfy the following recovery guarantees. Consider the same setting as in Proposition~\ref{prop:compressible_eta} and assume that the Riesz constants $b_{\Phi}$ and $B_{\Phi}$ in \eqref{eq:Riesz_constants_eta_compressible} satisfy the  sufficient condition
\begin{equation}
\label{eq:suff_cond_on_Riesz_constants}
\frac{b_\Phi}{B_\Phi} > 
\begin{cases}
1- \frac{0.98}{13} \approx 0.9246, & \text{(OMP)}\\
1-\frac{0.98}{\sqrt{2}}\approx 0.3070. & \text{(QCBP)}
\end{cases}
\end{equation}
Let $s\in \mathbb{N}$ and assume that the number of collocation points satisfies
\begin{equation}
\label{eq:suff_cond_m_CFC}
m \geq c \max\{B_\Phi^{-2},1\} K_\Phi^2 \, s \, L,
\end{equation}
where
\begin{equation}
\label{eq:log_factor_CFC}
L = \log^2(sK_\Phi^2\max\{B_\Phi^{-2},1\}) \min\{\log(n) + d, \log(2n) \log(2d)\} + \log(2\epsilon^{-1}),
\end{equation}
and $K_{\Phi}$ is defined as in \eqref{eq:unif_bound_eta_compressible}. Let $\hat{\bm{c}} \in \mathbb{C}^N$ be either the OMP solution computed via $K = 24s$ iterations or any QCBP solution with tuning parameter $\gamma \geq \|\bm{e}\|_2/\sqrt{B_{\Phi}}$ where $\bm{e} \in \mathbb{C}^m$ is as in \eqref{eq:def_error_vector}. Then, the corresponding compressive Fourier collocation approximation $\hat{u}$ to $u$ satisfies the following error bounds for all $f\in L^2(\mathbb{T}^d)$ with probability at least $1-\epsilon$:
\begin{align}
\label{eq:CFC_L2_Laplacian_bound}
\|\Delta(u-\hat{u}) \|_{L^2} 
& \leq \|\Delta(u-u_{\Lambda})\|_{L^2} +C_1\frac{ \sigma_s(\bm{c}_\Lambda)_1}{\sqrt{s}} + C_2 \|\nabla \cdot (a \nabla (u - u_{\Lambda}))\|_{L^\infty},\\
\label{eq:CFC_L2_bound}
\|u-\hat{u} \|_{L^2} 
& \leq \|u-u_{\Lambda}\|_{L^2} +  \frac{1}{4\pi^2} \left(C_1 \frac{ \sigma_s(\bm{c}_\Lambda)_1}{\sqrt{s}} + C_2 \|\nabla \cdot (a \nabla (u - u_{\Lambda}))\|_{L^\infty}
\right).
\end{align}
%Moreover, 
%$$
%\|\bm{e}\|_2 \leq \|\nabla \cdot (a \nabla (u - %u_{\Lambda}))\|_{L^\infty}. 
%$$
%\rev{The convergence theorem of OMP and QCBP looks similar. But there are difference in constants $C_1, C_2, C_3, C_4$.
%$$
%C_3, C_4 \sim Constant \cdot(1-\sqrt{2}\delta)
%$$
%}
\end{theorem}
%Some remarks
%\begin{itemize}
%    \item If the solution is sparse on the Fourier basis and $u \in \Span\{\Psi_\bm{\nu}\}$, then the truncation error equals to zero.
%\end{itemize}
The error bounds \eqref{eq:CFC_L2_Laplacian_bound}--\eqref{eq:CFC_L2_bound} show that the compressive Fourier collocation method is able to approximate compressible solutions in an \emph{accurate} and \emph{stable} way. Accuracy is due to the presence of the best $s$-term approximation error term $\sigma_2(\bm{c}_{\Lambda})/\sqrt{s}$ (optimal up to constants, due to the instance optimality theory; see, e.g., \citet[Chapter 11]{foucart2013mathematical}) and of the truncation errors involving the term $u-u_{\Lambda}$. Stability is intended to be with respect to the diffusion coefficient $a$ (defining the error term $\bm{e}$ in \eqref{eq:def_error_vector}) and is implied by the presence of $\|\nabla\cdot(a\nabla(u-u_{\Lambda}))\|_{L^\infty}$. 

Disregarding the dependence on the probability of failure $\epsilon$ and the constants $b_{\Phi}$, $B_{\Phi}$ and $K_{\Phi}$, the sufficient condition \eqref{eq:suff_cond_m_CFC} reads
$$
m \gtrsim s \log^2(2s) \min\{\log(n) + d, \log(2n) \log(2d)\}.
$$
Hence, the dependence of the number of collocation points $m$ on the dimension $d$ is logarithmic. This justifies our claim that the compressive Fourier collocation method is able to lessen the curse of dimensionality with respect to the number of collocation points under suitable regularity conditions on the diffusion coefficient $a$. %Let us add some comments on these conditions. They are very conservative. The numerics will show that in practice they do not seem to be so restrictive, although this issue deserves further investigation.
\begin{remark}[Removing condition \eqref{eq:suff_cond_on_Riesz_constants} for QCBP]
Condition \eqref{eq:suff_cond_on_Riesz_constants} can be removed in the QCBP case at the cost of multiplying $C_2$ by $B_{\Phi}/b_{\Phi}$ in the error bounds using the robust null space property approach in \cite[Theorem 2.6]{brugiapaglia2021sparse}. Applying this result, the condition on the number of collocation points becomes
$$
m \geq c \left(\frac{\max\{1,B_{\Phi}\}}{b_{\Phi}}\right)^2 K_{\Phi}^2\, s\,  L,
$$
where
$$
L = \log^2\left(sK_\Phi^2\frac{\max\{1,B_\Phi^{2}\}}{b_{\Phi}}\right) \min\{\log(n) + d, \log(2n) \log(2d)\} + \log(2\epsilon^{-1}),
$$
and the constant $C_2$ gets multiplied by $(B_{\Phi}/b_{\Phi})$. Therefore, condition \eqref{eq:suff_cond_on_Riesz_constants} does not seem to be a fundamental limit on the applicability of the compressive Fourier collocation methods and its ability to lessen the curse of dimensionality in the number of collocation points.
\end{remark}

Finally, we note that the sufficient condition on the number of collocation points in Theorem~\ref{thm:recovery_CFC} is significantly weaker with respect to the corresponding condition in \citet[Theorems 3 and 4]{brugiapaglia2020compressivespectral}, where there is an explicit dependence of $m$ on $2^d$.

\section{Numerical  experiments}
\label{s:numerics}

In this section, we present numerical experiments on the the model problem \eqref{eq:diffusion_eq_periodic_1}--\eqref{eq:diffusion_eq_periodic_2} using the compressive Fourier collocation method presented in Section~\ref{s:CFC_method}. In \S\ref{s:sol_coeff}, we start by defining the exact solutions (sparse and nonsparse) and the diffusion coefficients (constant, sparse, and nonsparse) considered in the experiments. Then, we numerically demonstrate the accuracy and stability of the method (already  established theoretically by Theorem~\ref{thm:recovery_CFC}) in dimension $d = 2$ (\S\ref{s:2D}) and $d=8$ and $d=20$ (\S\ref{s:8D}). We also study the probability of success as a function of the number of collocation points for sparse solutions in Appendix~\ref{app:phase_transition}. The main aspects of our experimental setting are described in this section. However, further technical details are provided in Appendix~\ref{a:details_numerics}. %We provide further details on the numerical experiment%We also represent an example with different index sets to illustrate that the size of the index sets is a crucial factor to the error. 

\paragraph{Measurement of errors.} In all numerical experiments, we use a relative $L^2$-error to measure the approximation error. It is defined as
$$
\text{relative}\ L^2\text{-error} = \frac{\|u-\hat{u}\|_{L^2}}{\|u\|_{L^2}},
$$
where $\hat{u}$ is the computed approximation to an exact solution $u$. Considering that the $L^2$-norms are unapproachable in the high-dimensional case, $\|\cdot\|_{L^2}$ are evaluated using Monte Carlo method, i.e.
$$
\|u\|_{L^2} \approx \sqrt{\frac{1}{M}\sum^{M}_{i=1} |u(\bm{x}_i)|^2},
$$
where $\bm{x}_i$ are $M \gg 1$ points randomly distributed in ${\mathbb{T}^d}$. We choose $M=2|\Lambda|$ in practice, with the exception of the last experiment in \S\ref{app:phase_transition} (Probability of successful recovery) where $M=200$.

\subsection{Exact solutions and diffusion coefficients} 
\label{s:sol_coeff}

We consider the following exact solutions and diffusion coefficients. For now, we assume $d = 2$.

\begin{figure}[!t]
\begin{subfigure}{.48\textwidth}
    \centering
    \includegraphics[width=\linewidth]{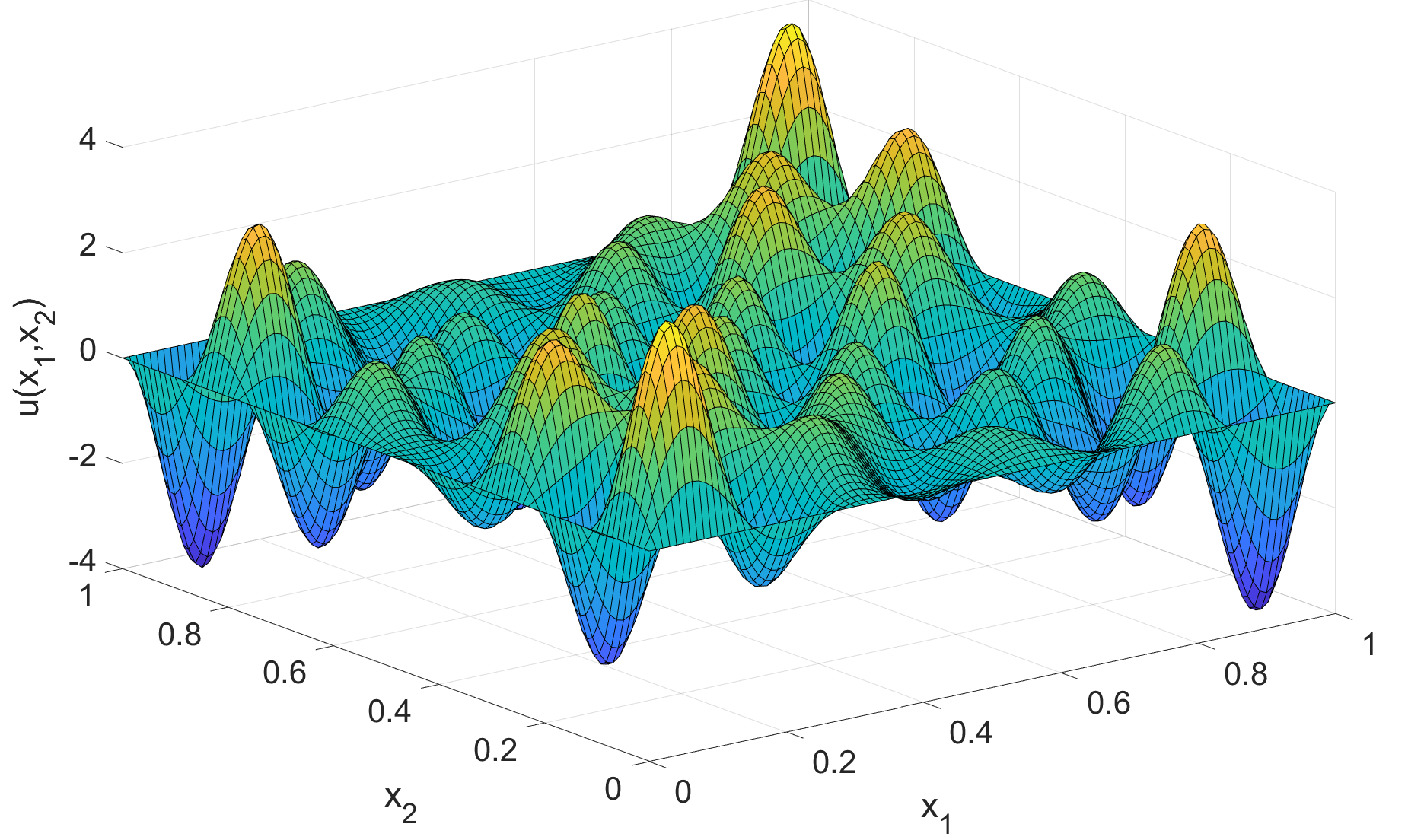}
    \caption{$u_1$}
\end{subfigure}
\begin{subfigure}{.48\textwidth}
    \centering
    \includegraphics[width=\linewidth]{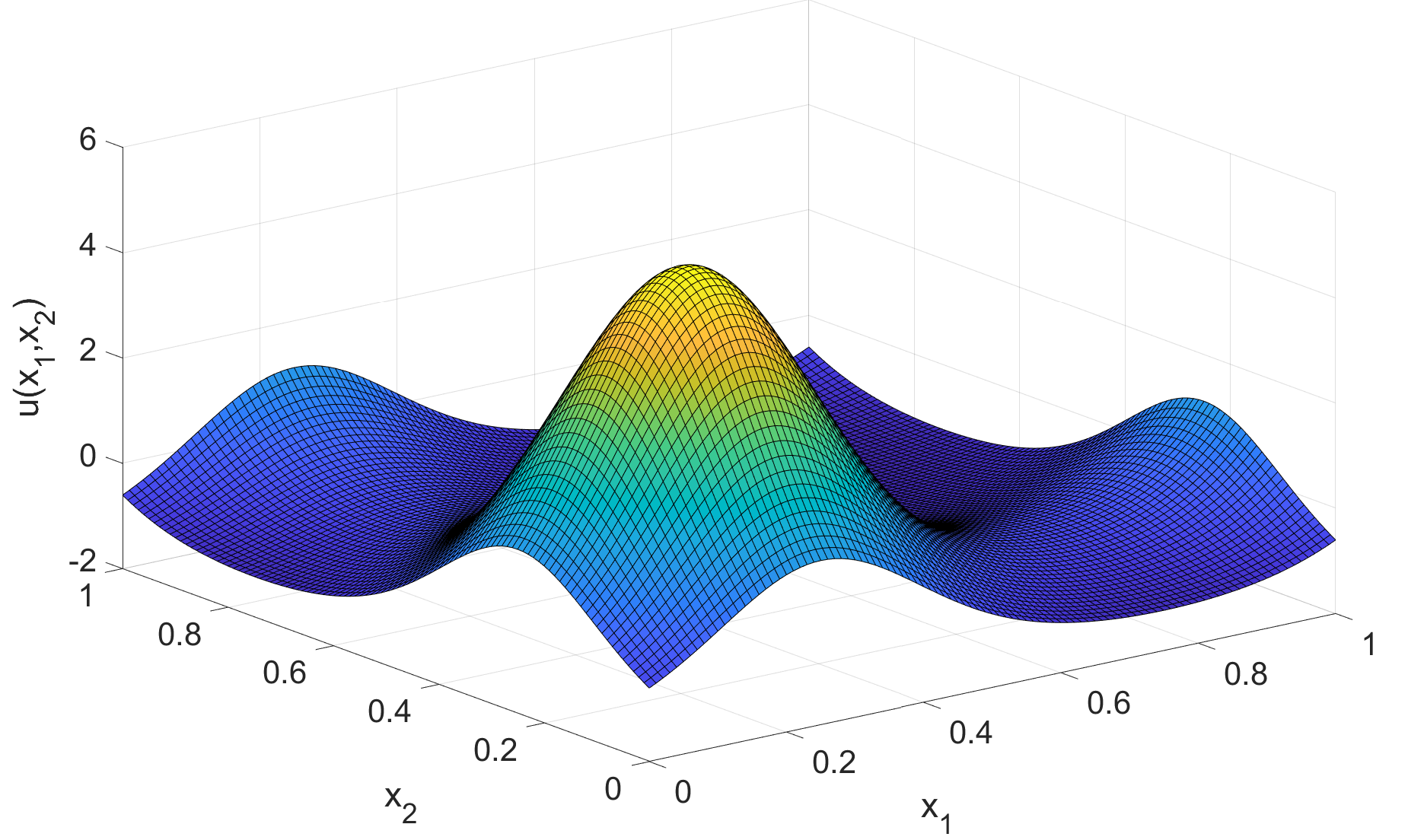}
    \caption{$u_2$}
\end{subfigure}
\caption{\label{fig:exact}The sparse and nonsparse exact solutions $u_{1}$ and $u_{2}$, defined in  \eqref{eq:exact_sparse_solu} and \eqref{eq:exact_nonsparse_solu}, respectively. }
\end{figure}

\paragraph{Exact solutions.}
We consider both sparse and nonsparse exact solutions. The sparse solutions (with respect to the Fourier basis) are randomly generated as 
\begin{equation} \label{eq:exact_sparse_solu}
u_{1}(\bm{x}) = \sum^{10}_{k=1} d_k \sin{(2\pi m_k x_1)}\sin{(2\pi n_k x_2)}, \quad \forall \bm{x} \in \mathbb{T}^d,
\end{equation} 
where $d_i$ are random independent real coefficients uniformly distributed in $[0,1]$, and $m_k$, $n_k$ are random integers such that $\sin{(2\pi m_k x_1)}\sin{(2\pi n_k x_2)} \in \Span\{F_{\bm{\nu}}\}_{\bm{\nu}\in \Lambda}$ (recall that $\Lambda$ is the truncation set used in the compressive Fourier collocation scheme; see also Algorithm~\ref{alg:Compressive_Fourier_collocation} ). The function $u_1$ is $40$-sparse with respect to the Fourier basis. We also consider a nonsparse solution defined as 
\begin{equation} \label{eq:exact_nonsparse_solu}
u_{2}(\bm{x}) = \exp{\left(\sin{(2\pi x_1)}+\sin{(2\pi x_2)}\right)} - c, \quad \forall \bm{x} \in \mathbb{T}^d.
\end{equation}
Here, the constant $c$ is chosen in such a way to guarantee that $\int_{\mathbb{T}^d} u_{2}(\bm{x}) \mathrm{d} \bm{x}=0$, i.e.\ to satisfy condition \eqref{eq:diffusion_eq_periodic_2} (note that this is always the case for $u_{1}$). These exact solutions are shown in Fig.~\ref{fig:exact}. 

\paragraph{Diffusion coefficients.}
In all experiments, we consider three types of diffusion coefficients: namely, constant, sparse, and non-sparse diffusion coefficients. They are defined as follows, for any $\bm{x} \in \mathbb{T}^d$:
\begin{align}
    a_1(\bm{x}) & =1, \\
    a_2(\bm{x}) & = 1+ 0.25\sin{(2\pi x_1)}\sin{(2\pi x_2)} + 0.25\sin{(4\pi x_1)}, \\
    a_3(\bm{x}) & =1 + 0.2 \exp{(\sin{(2\pi x_1)} \sin{(2\pi x_2)})}.
\end{align}  
We are interested in assessing the stability of the method over different types of diffusion.

\subsection{Low-dimensional experiments ($d=2$)}
\label{s:2D}
For the two-dimensional experiments, we consider the hyperbolic cross index set $\Lambda$ (minus the zero index) defined as in \eqref{eq:def_Lambda}, with $d=2$ and of order $n=39$. This leads to $N = |\Lambda| = 445$. %$s$ is the number of the iterations of the OMP method, and $m$ is the number of collocation points (using Monte Carlo) for all methods. We tune the parameter in QCBP using the cross validation,  which is generally the optimized choice. \rev{QCBP parameter reference here.} 

\begin{figure}[!t]
\begin{subfigure}{.36\textwidth}
  \centering
  \includegraphics[height=38mm]{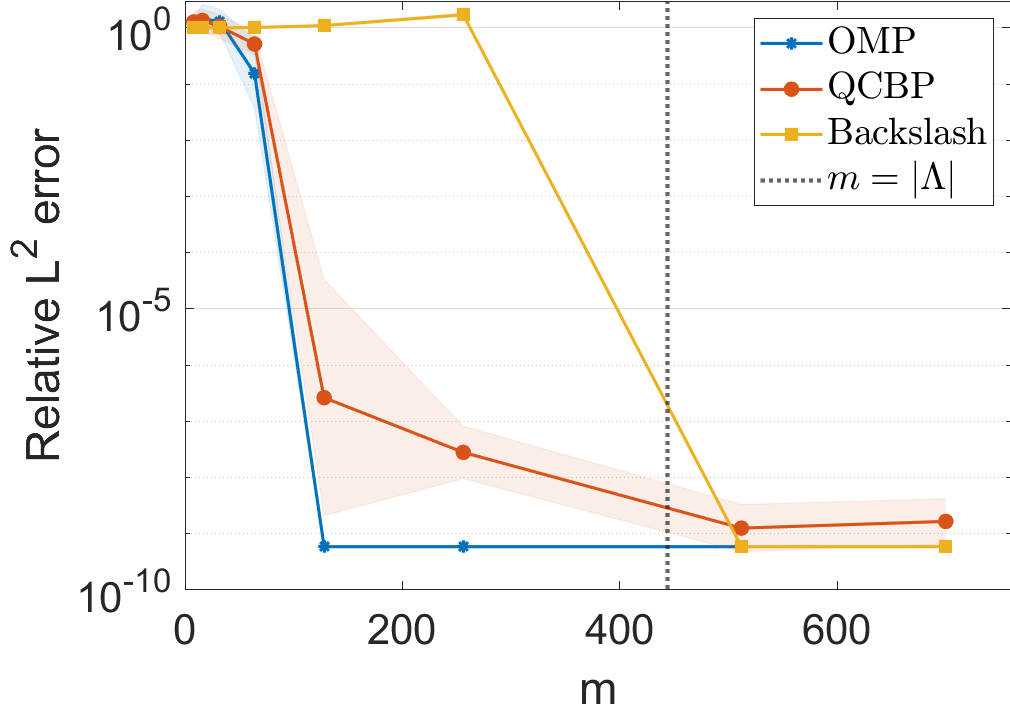}
  \caption{constant diffusion $a_1$}
\end{subfigure}
\begin{subfigure}{.30\textwidth}
  \centering
  \includegraphics[height=38mm]{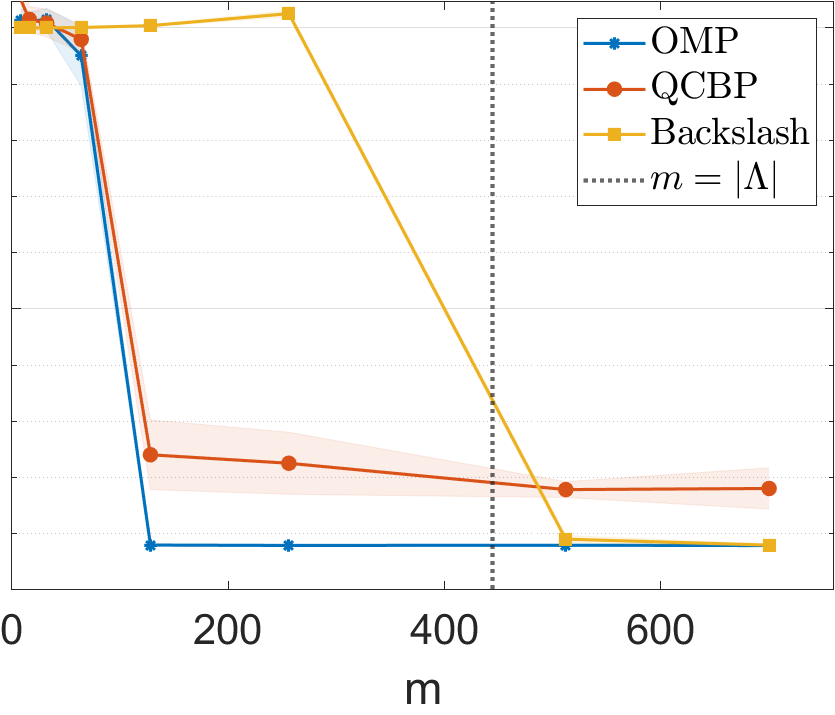}
  \caption{sparse diffusion $a_2$}
\end{subfigure}
\begin{subfigure}{.30\textwidth}
  \centering
  \includegraphics[height=38mm]{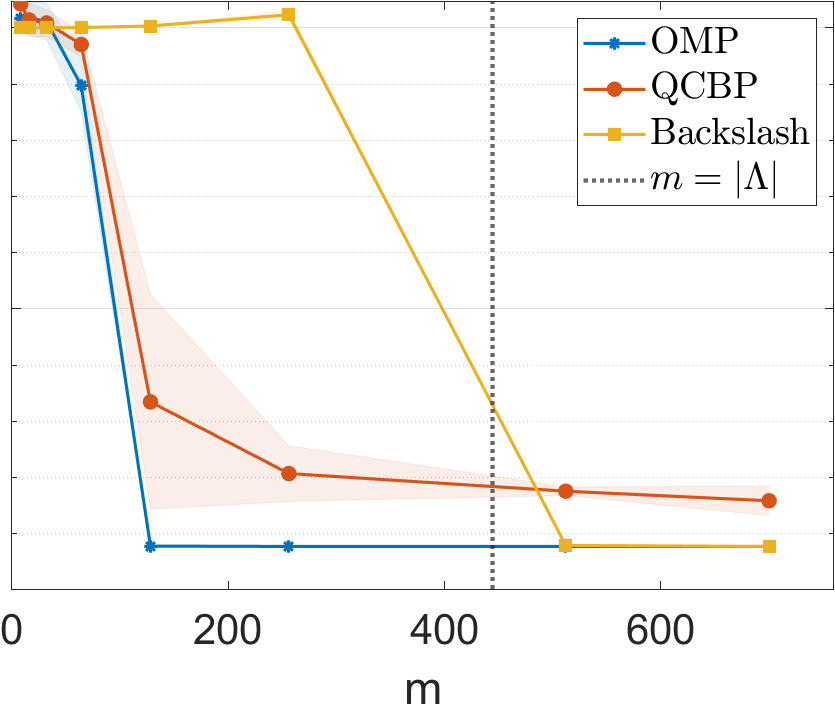}
  \caption{nonsparse diffusion $a_3$}
\end{subfigure}
\caption{($d=2$, sparse solution) Relative $L^2$-error versus number of collocation points $m$ for different diffusion coefficients. The exact solution is $u_{1}$, defined in \eqref{eq:exact_sparse_solu}.} \label{Fig:2D_sparse}
\end{figure}

% Sparse exact solution and nonsparse exact solution are corresponding to $u_{\mbox{exact}1}$ and $u_{\mbox{exact}2}$.
\paragraph{Sparse solution ($d=2$).}
We test compressive Fourier collocation with the sparse exact solution $u_{1}$ defined in \eqref{eq:exact_sparse_solu} and where recovery is performed via the (MATLAB\textregistered) backslash operator (corresponding to ordinary least squares \eqref{eq:OLS}), OMP (Algorithm~\ref{alg:OMP}), and QCBP (defined in \eqref{eq:QCBP}). The frequencies $m_i$ and $n_i$ are chosen as independent random integers uniformly distributed in $\{1,\ldots,5\}$. This choice ensures, on the one hand, that $\sin{(2\pi m_i x_1)}\sin{(2\pi n_i x_2)} \in \Span\{F_{\bm{\nu}}\}_{\bm{\nu} \in \Lambda}$ and, on the other hand, that the error arising from the computation of  the right-hand side $\bm{b}$ via finite differences is negligible (see Appendix~\ref{a:details_numerics} for further details). Recall that $u_1$ is 40-sparse with respect to the complex Fourier basis. Fig.~\ref{Fig:2D_sparse} shows the relative $L^2$-error associated with OMP, QCBP, and backslash as a function of the number of collocation points $m$. 
OMP and QCBP recover the exact solution with error less than $10^{-6}$ for $m=128 \ll |\Lambda|$, which is only a small multiple of the target sparsity $s = 40$. The dashed vertical line indicates $m=|\Lambda|$, where the number of collocation points is the same as the number of unknowns. In order to produce accurate results, the backslash (least squares method) needs at least $m\geq |\Lambda|$ collocation points. %For $m\geq |\Lambda|$, the linear system becomes overdetermined, and the backslash method (least squares) also produces accurate results. 
Results in Fig.~\ref{Fig:2D_sparse} illustrate that the Compressive Fourier collocation methods considered are able to accurately recover sparse exact solutions to the periodic diffusion equation and are stable with respect to variations of the diffusion coefficient. 

\begin{figure}[t]
\begin{subfigure}{.36\textwidth}
  \centering
  \includegraphics[height=38mm]{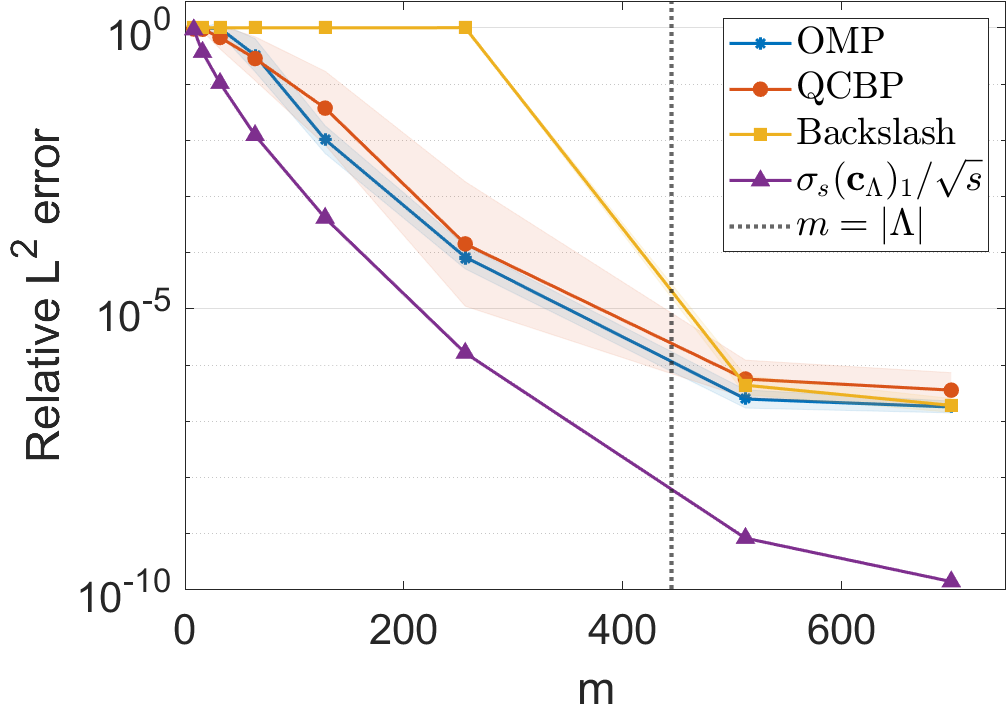}
  \caption{constant diffusion $a_1$}
\end{subfigure}
\begin{subfigure}{.30\textwidth}
  \centering
  \includegraphics[height=38mm]{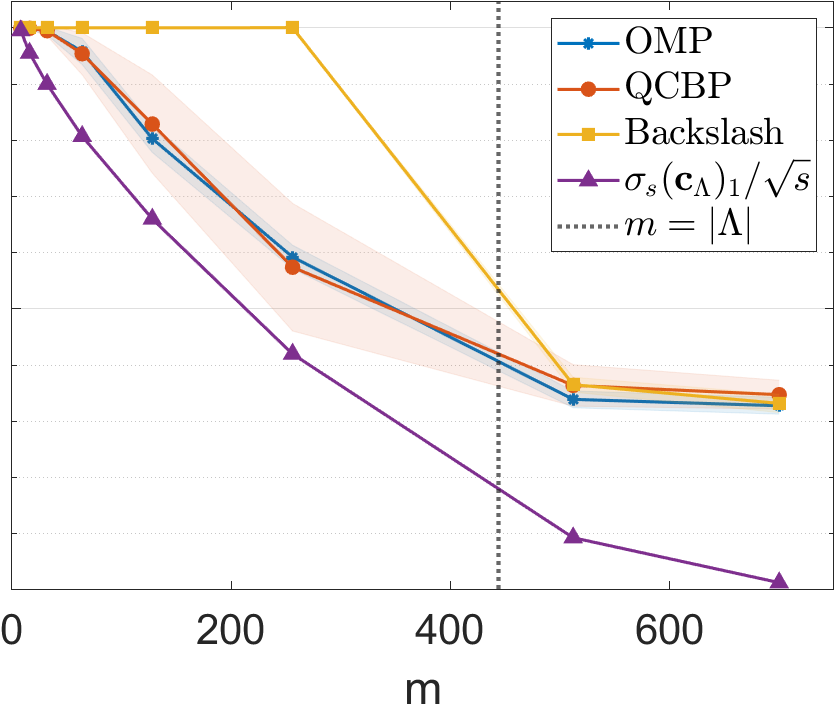}
  \caption{sparse diffusion $a_2$}
\end{subfigure}
\begin{subfigure}{.30\textwidth}
  \centering
  \includegraphics[height=38mm]{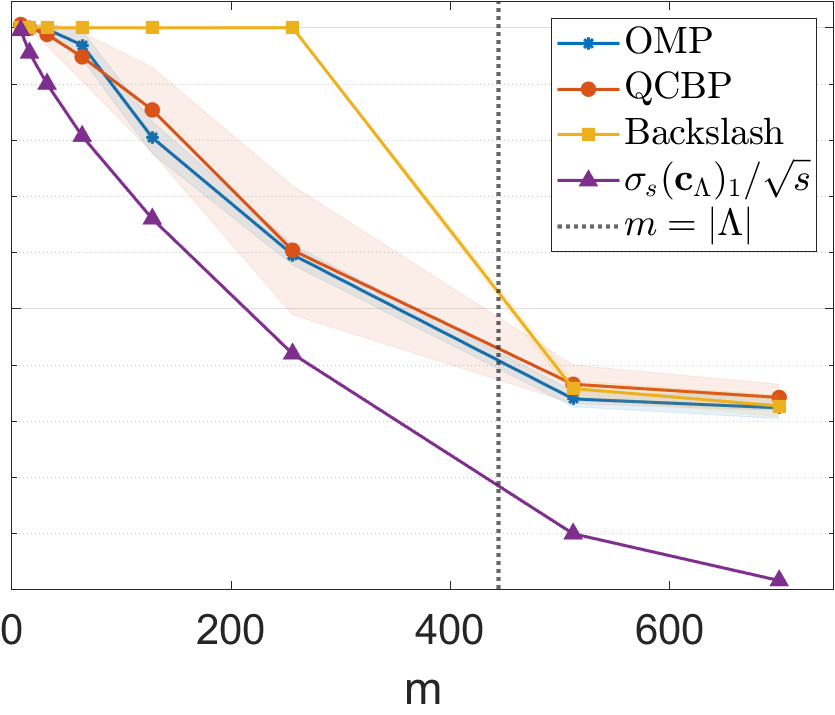}
  \caption{nonsparse diffusion $a_3$}
\end{subfigure}
\caption{($d=2$, nonsparse solution) Relative $L^2$-error versus number of collocation points $m$. The exact solution is $u_{2}$, defined in \eqref{eq:exact_nonsparse_solu}. \label{Fig:2D_nonsparse}}
\end{figure}

\paragraph{Nonsparse solution ($d=2$).}
In this test, the exact solution is the nonsparse solution $u_{2}$ defined in \eqref{eq:exact_nonsparse_solu}. The parameter settings are the same as in the sparse solution test. The results are shown in Fig.~\ref{Fig:2D_nonsparse}.
The purple lines represent the term $\sigma_s(\bm{c}_\Lambda)_1/\sqrt{s}$ (that appears in error bound of Theorem~\ref{thm:recovery_CFC}) with $s = m/2$ and where $\bm{c}_\Lambda$ is approximated using compressive Fourier collocation with $m = 4|\Lambda|$ via the backslash. The fast decay of $\sigma_s(\bm{c}_\Lambda)_1/\sqrt{s}$ shows that $u_2$ is nonsparse. Moreover, the results of the OMP method with $K = s$ iterations are at most $s$-sparse. Thus, the purple curve is a lower bound for the error of the OMP method (thanks to Theorem~\ref{thm:RIP->recovery}). The results of Fig.~\ref{Fig:2D_nonsparse} demonstrate that the compressive Fourier collocation solutions recovered via OMP and QCBP approach the exact solution as the number of collocation points $m$ increases, irrespective of the sparsity of the diffusion coefficient. 

\subsection{High-dimensional experiments ($d = 8$ and $d = 20$)}
\label{s:8D}
In the next experiments, we consider a high-dimensional setting with $d=8$ and $d=20$. The index set $\Lambda$ naturally differs from that of the case $d=2$. The exact solution is also slightly different from the one considered in the case $d=2$. The other  parameters are as in \S\ref{s:2D}.

\begin{figure}[!t]
\begin{subfigure}{.36\textwidth}
  \centering
  \includegraphics[height=38mm]{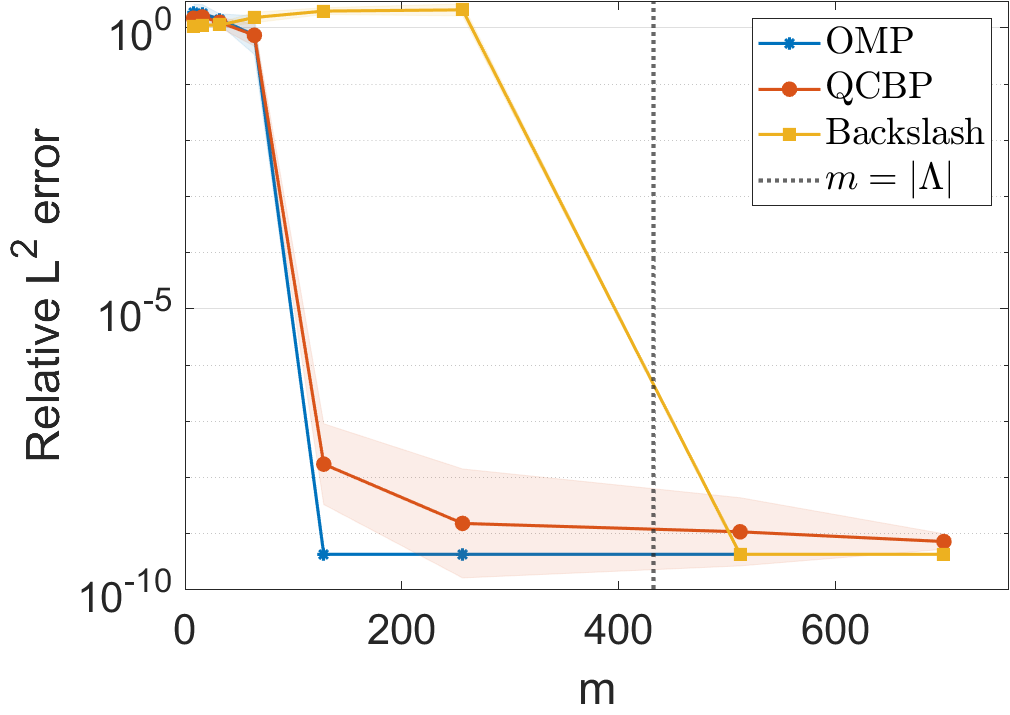}
  \caption{constant diffusion $a_1$}
\end{subfigure}
\begin{subfigure}{.30\textwidth}
  \centering
  \includegraphics[height=38mm]{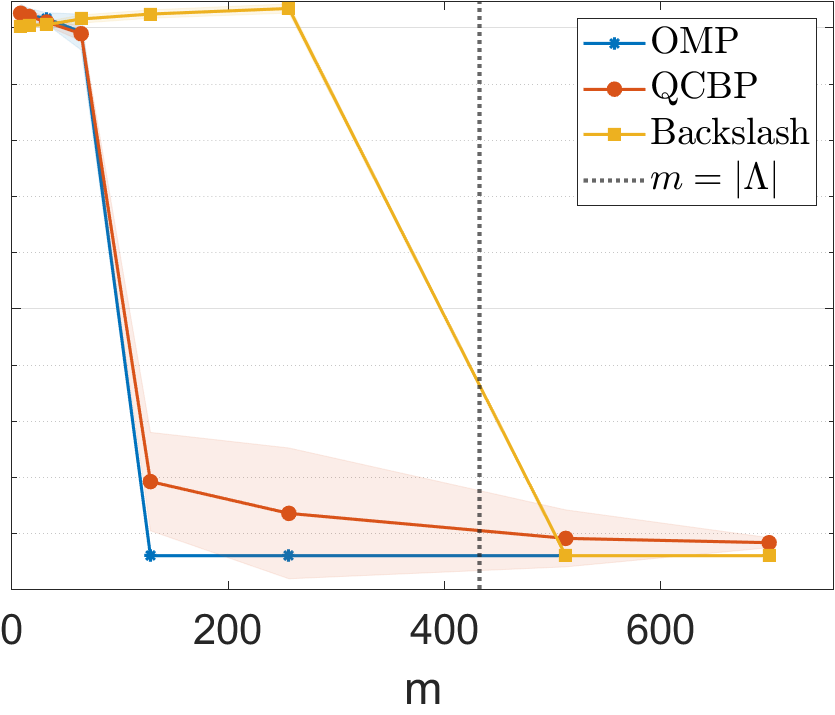}
  \caption{sparse diffusion $a_2$}
\end{subfigure}
\begin{subfigure}{.30\textwidth}
  \centering
  \includegraphics[height=38mm]{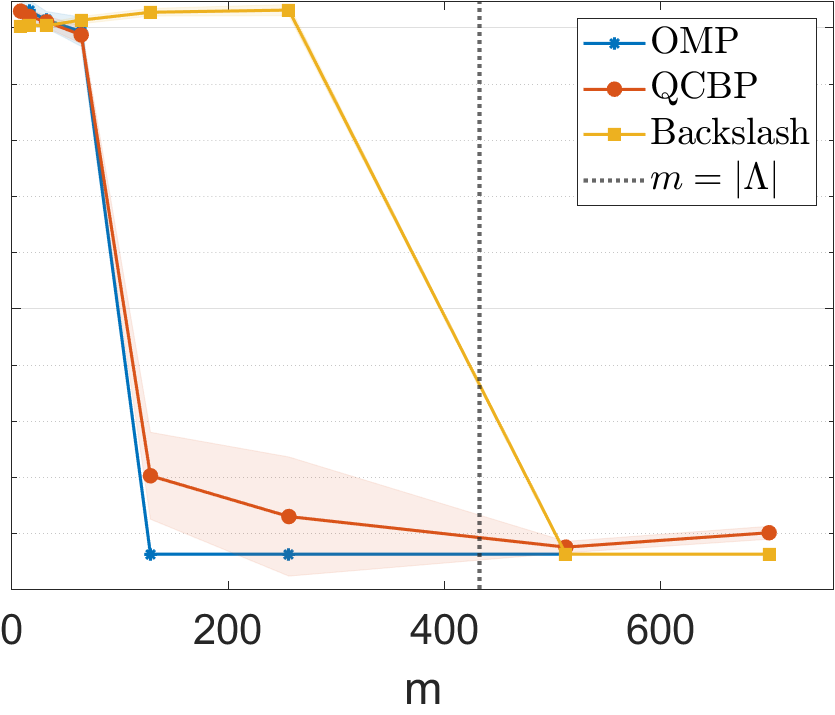}
  \caption{nonsparse diffusion $a_3$}
\end{subfigure}
\caption{($d=8$, sparse solution) Relative $L^2$-error versus number of collocation points $m$. The exact solution is $u_{3}$, defined in \eqref{eq:def_u3}.} \label{Fig:8D_sparse}
\end{figure}

\paragraph{Sparse exact solution ($d=8$).}
We choose the index set $\Lambda$ as the $8$-dimensional hyperbolic cross in \eqref{eq:def_Lambda} of order $n=7$. This leads to $N = |\Lambda|=432$. We consider a sparse exact solution given by
\begin{equation}
\label{eq:def_u3}
u_{3}(\bm{x}) = \sum^{10}_{k=1} d_k \prod_{l=1}^8 \sin(2\pi m_{k,l} x_l), \quad \forall\bm{x} \in \mathbb{T}^8,
\end{equation}
%where $s=10$ is the sparsity with respe. The sparsity of the solution on the complex Fourier basis is less than 40. 
where the frequencies $\bm{m}_k$ are randomly and uniformly selected from the index set $\Lambda$ and $d_i$ are random real coefficients uniformly distributed in $[0,1]$. Note that since $n=7$ we have $\|\bm{m}_k\|_0 \leq 2$. Hence, $u_3$ is 40-sparse with respect to the complex Fourier system. Other parameters are the same as in the case $d=2$ for the sparse solution test.  Fig.~\ref{Fig:8D_sparse} shows that the accuracy of compressive Fourier collocation is independent of the number of dimensions when the exact solution is sparse. %The OMP method can achieve the desired sparsity with $m \approx 2s$.  

\begin{figure}[!t] 
\begin{subfigure}{.36\textwidth}
  \centering
  \includegraphics[height=38mm]{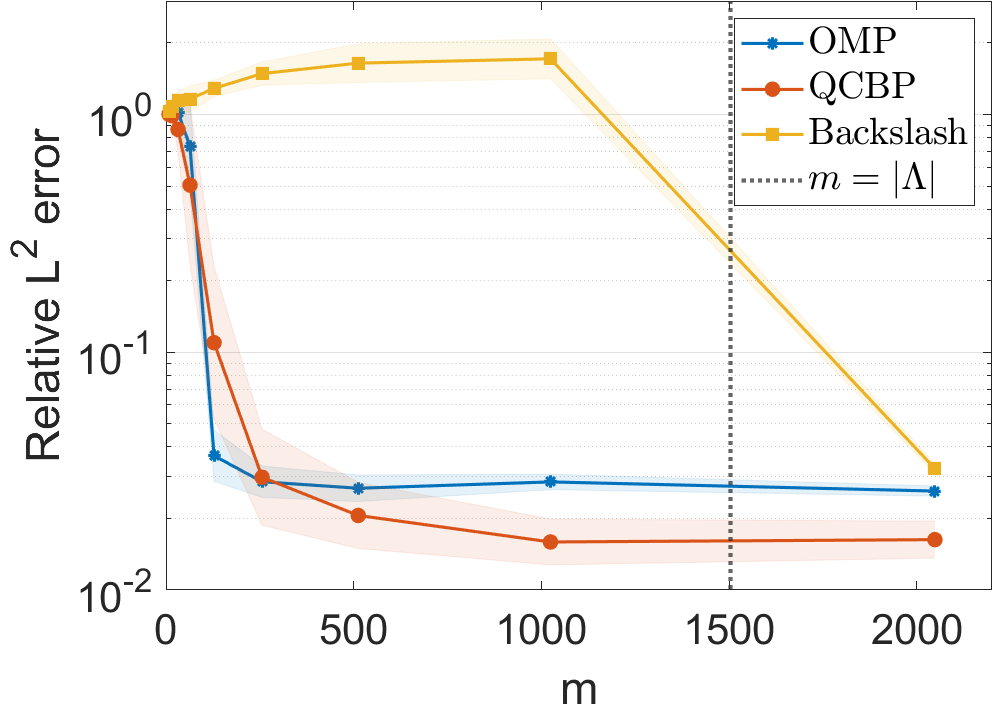}
  \caption{constant diffusion $a_1$}
\end{subfigure}
\begin{subfigure}{.30\textwidth}
  \centering
  \includegraphics[height=38mm]{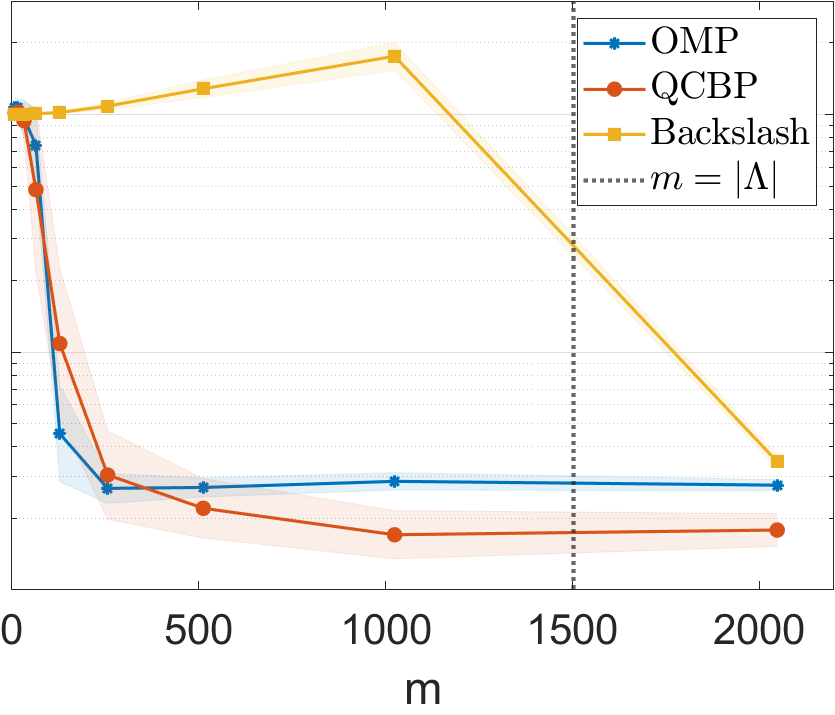}
  \caption{sparse diffusion $a_2$}
\end{subfigure}
\begin{subfigure}{.30\textwidth}
  \centering
  \includegraphics[height=38mm]{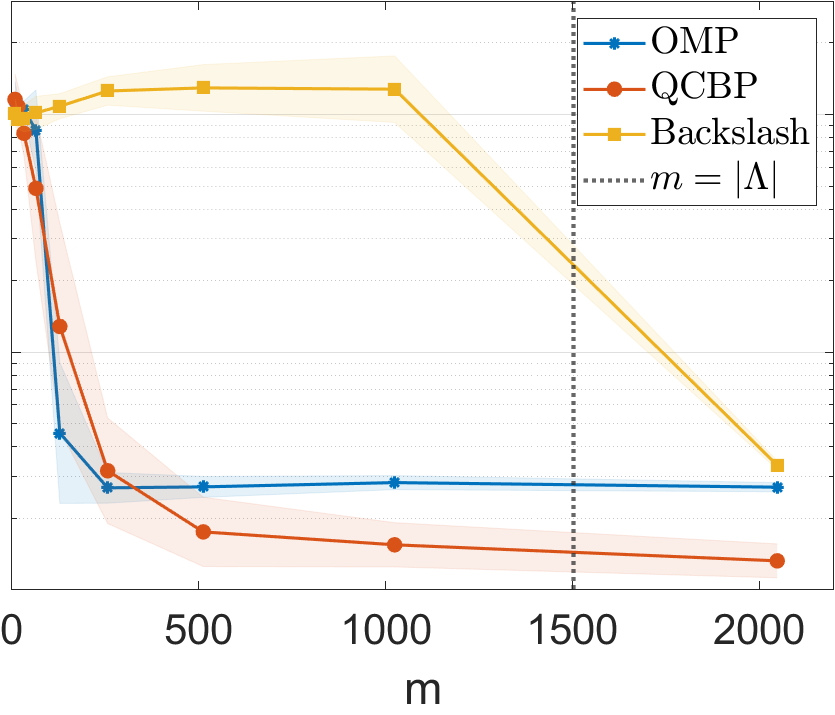}
  \caption{nonsparse diffusion $a_3$}
\end{subfigure}
\caption{($d=8$, nonsparse solution) Relative $L^2$-error versus number of collocation points $m$. The exact solution is $u_{2}$, defined in \eqref{eq:exact_nonsparse_solu}.} \label{Fig:8D_nonsparse}
\end{figure}

\paragraph{Nonsparse solution ($d=8$).}
We now consider the exact solution $u_{2}$ defined in \eqref{eq:exact_nonsparse_solu} but in dimension  $d=8$ (only the first two variables are active). We let $\Lambda$ be the 8-dimensional hyperbolic cross in \eqref{eq:def_Lambda} of order $n=11$ (corresponding to $N = |\Lambda|=1505$).  As shown in Fig.~\ref{Fig:8D_nonsparse}, the errors for all three numerical methods are significantly larger than the error in the two-dimensional case due to the high dimensionality. However, QCBP and OMP are able to approximate the solution with a relative error close to $10^{-2}$ using only a small amount of collocation points and regardless of the sparsity of the diffusion coefficients.

\paragraph{The impact of the index set $\Lambda$ ($d=8$).} We consider the same settings as in the nonsparse solution case with $d=8$ and the nonsparse diffusion coefficient $a_3$ and consider index sets $\Lambda$ of increasing size. We use this example to illustrate the impact of the size of the index set. Since the error of the OMP and QCBP reach a plateau for $m = 1024 \ll |\Lambda|$ in Fig.~\ref{Fig:8D_nonsparse}, we only present the results with $m \leq 1024$ to reduce the computational cost of the experiment. 
\begin{figure}[!t]
\begin{subfigure}{.36\textwidth}
  \centering
  \includegraphics[height=38mm]{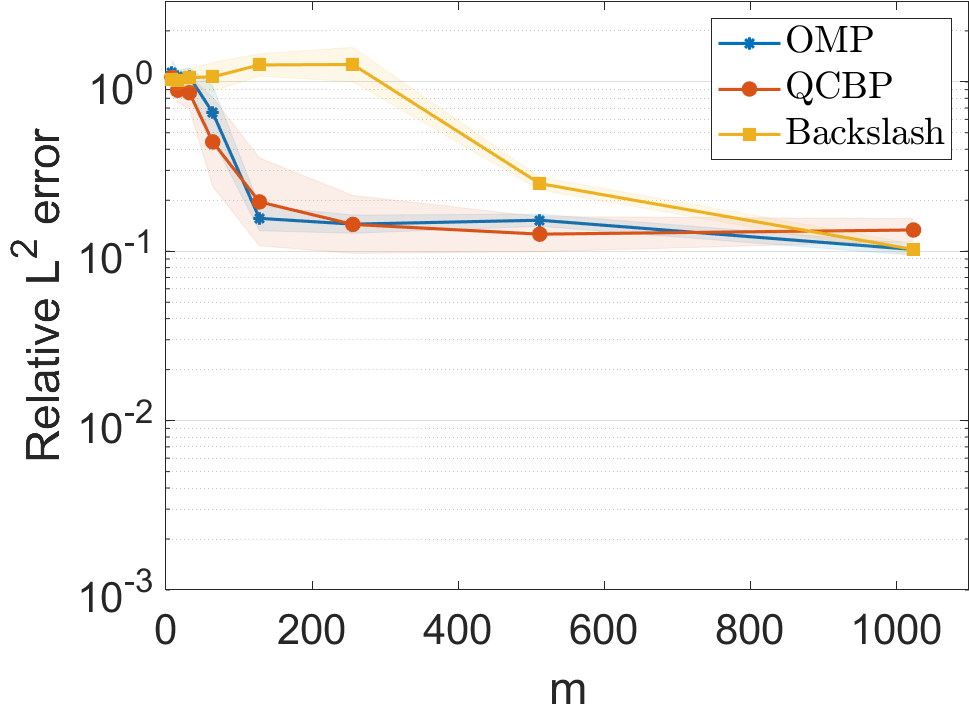}
  \caption{$n=7$}
\end{subfigure}
\begin{subfigure}{.30\textwidth}
  \centering
  \includegraphics[height=38mm]{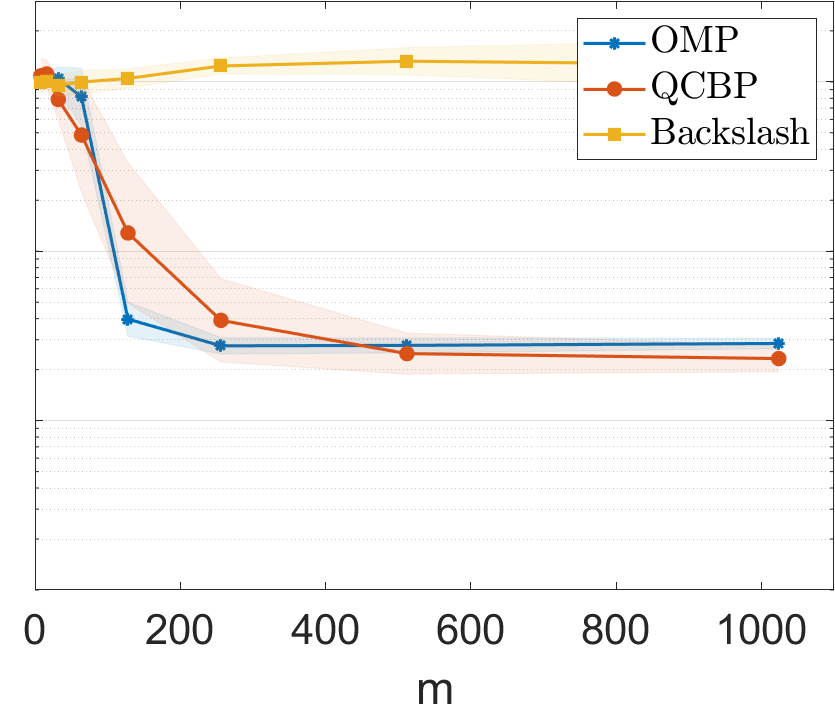}
  \caption{$n=11$}
\end{subfigure}
\begin{subfigure}{.30\textwidth}
  \centering
  \includegraphics[height=38mm]{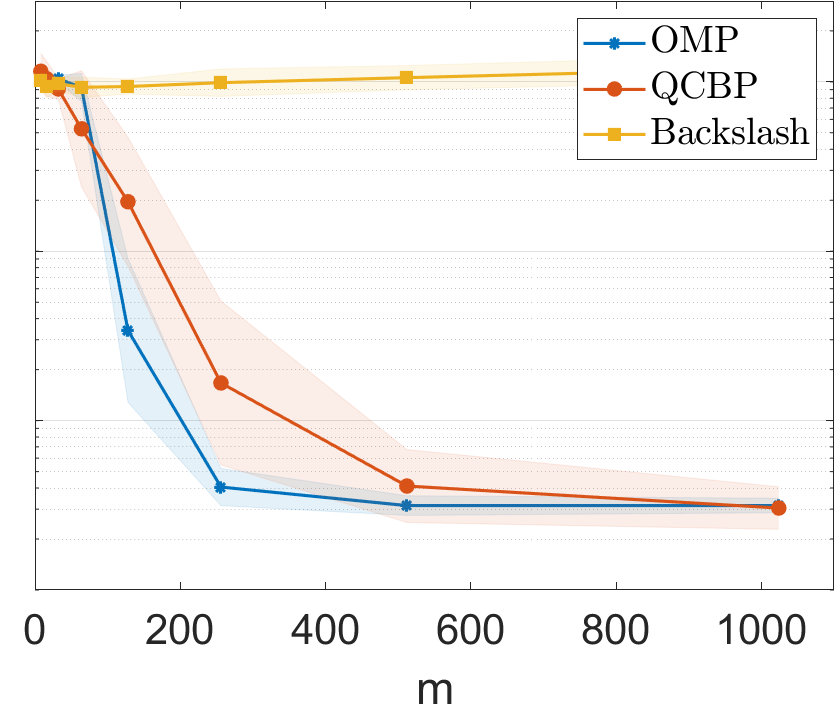}
  \caption{$n=15$}
\end{subfigure}
\caption{($d=8$, impact of $\Lambda$) Relative $L^2$-error versus number of collocation points $m$. The exact solution, $u_2$, and the diffusion coefficient, $a_3$, are nonsparse. The cardinality of the index set $N = |\Lambda|$ is $432$ in (a),  $1505$ in (b), and $3808$ in (c).} \label{Fig:8D_indexset}
\end{figure}

We compare the performance of the compressive Fourier collocation method for index sets $\Lambda$ of different sizes. Fig.~\ref{Fig:8D_indexset} shows that as the size of the index set increase, the error of OMP and QCBP is reduced from around $10^{-1}$ to around $10^{-3}$. This shows that size of the index set is a key factor in reducing the approximation error in the high-dimensional setting.

\begin{figure}[!t]
\begin{subfigure}{.49\textwidth}
  \centering
  \includegraphics[height=50mm]{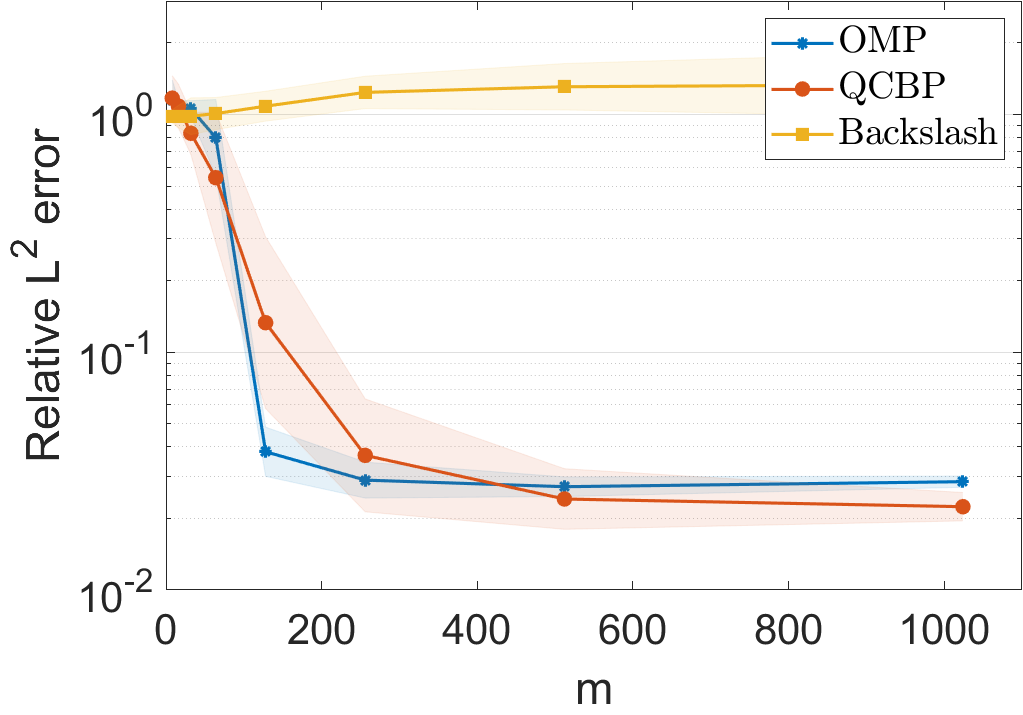}
  \caption{$d=8$}
\end{subfigure}
\begin{subfigure}{.49\textwidth}
  \centering
  \includegraphics[height=50mm]{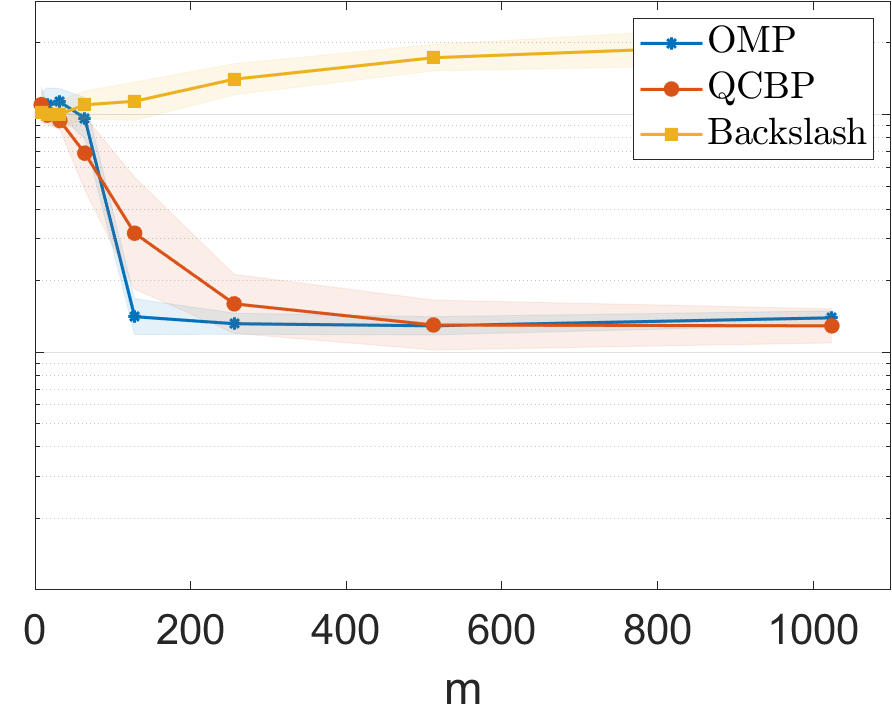}
  \caption{$d=20$}
\end{subfigure}
\caption{(The impact of dimensionality) Relative $L^2$-error versus number of collocation points $m$. The exact solution, $u_2$, and the diffusion coefficient, $a_3$, are nonsparse. The multi-index set $\Lambda$ is the largest hyperbolic cross set such that $|\Lambda| < 2550$ in (a) $d=8$ and (b) $d=20$.} \label{Fig:20D_indexset}
\end{figure}

\paragraph{The impact of dimensionality ($d=8$ and $d=20$).} %We now test compressive Fourier collocation on a higher-dimensional example. 
We now consider the same settings as in the nonsparse solution case and the nonsparse diffusion coefficient $a_3$, comparing the cases of dimension $d=8$ and $d=20$. To choose $\Lambda$, we consider the largest order $n$ such that the corresponding hyperbolic cross has cardinality less than a given budget value of 2550. This corresponds to $n=11$ (for $d=8$) and $n = 7$ (when $d=20$). The results are shown in  Fig.~\ref{Fig:20D_indexset}. We show that in dimension $d=20$, the approximation error increases and reaches a plateau at $10^{-1}$ because the multi-index set $\Lambda$ is only able to capture fewer significant coefficients due to the lower value of $n$ combined with the fact that the exact solution is highly anisotropic (it only depends on two variables out of 20).  Nonetheless, Fig.~\ref{Fig:20D_indexset} shows that the numerical solution approaches the exact solution in $d=20$, up to the truncation error corresponding to the choice of $\Lambda$. 

\paragraph{Probability of successful recovery versus the number of collocation points $m$.}
\label{app:phase_transition}
In this numerical experiment, we study the probability of successful sparse recovery with respect to the number of the collocation points $m$ for $d=2,8$. Similarly to equation~\eqref{eq:exact_sparse_solu}, the sparse solutions are defined as
\begin{equation} 
u(\bm{x}) = \sum^{q}_{k=1} d_k \sin{(2\pi m_k x_1)}\sin{(2\pi n_k x_2)}, \quad \forall \bm{x} \in \mathbb{T}^d,
\end{equation} 
where $d_k$ are random independent real coefficients uniformly distributed in $[0,1]$, $q$ is the sparsity respect to the real-valued Fourier basis, and $(m_k$, $n_k)$ are $q$ distinct random integer pairs in $\{1,2,3,4\}^2$ such that $\sin{(2\pi m_k x_1)}\sin{(2\pi n_k x_2)} \in \Span\{F_{\bm{\nu}}\}_{\bm{\nu}\in \Lambda}$.

\begin{figure}[!t]
\begin{subfigure}{.49\textwidth}
  \centering
  \includegraphics[height=50mm]{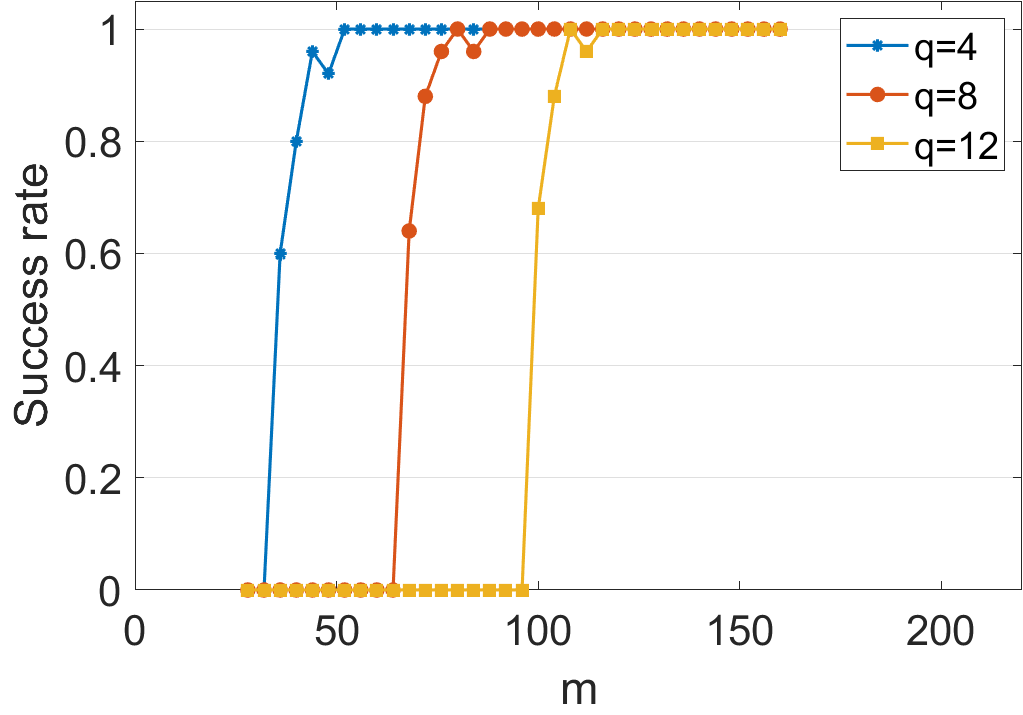}
  \caption{$d=2$}
\end{subfigure}
\begin{subfigure}{.49\textwidth}
  \centering
  \includegraphics[height=50mm]{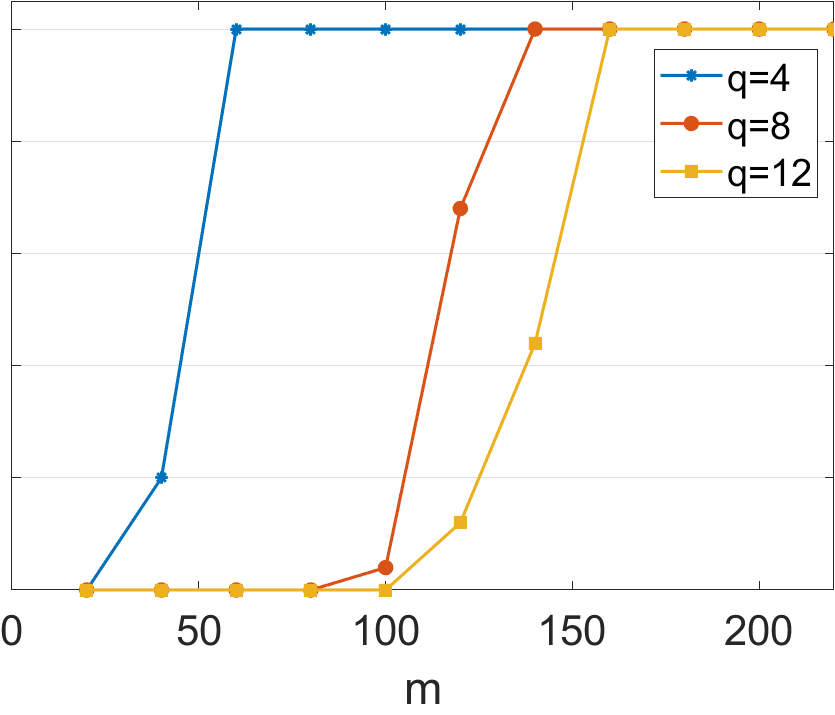}
  \caption{$d=8$}
\end{subfigure}
\caption{(Probability of successful recovery) Success rate of recovery versus number of collocation points $m$ for sparse solutions. The diffusion coefficient, $a_3$, is nonsparse.} \label{Fig:P_successful_recovery}
\end{figure}

We let $\Lambda$ be the hyperbolic cross set of order $n=26$ for both cases $d=2$ and $d=8$. The number of the collocation points is $m=2s$, where $s$ is the target sparsity and recovery is performed via the OMP algorithm (Algorithm \ref{alg:OMP} in Section~\ref{s:CFC_method}) with $K=s$ iterations. For each data point in Fig.~\ref{Fig:P_successful_recovery}, we perform 25 runs and compute the successful recovery rate by counting the number of runs corresponding to a relative $L^2$-error below $10^{-6}$. We compare the recovery success rate for sparse solutions with sparsity $q=4,\ 8,\ 12$. These sparse solutions have sparsity 16, 32, and 48 (i.e., $4q$) with respect to the complex Fourier basis, respectively. In order to reduce the computational cost of the experiment, we use a coarser grid in $m$ for $d=8$.

As shown in Fig.~\ref{Fig:P_successful_recovery}, the phase transition from unsuccessful to successful recovery happens at a value of $m$ slightly larger than $2s=8q$ for the two-dimensional case. In dimension $d=2$, the phase transition is quite sharp. However, in $d=8$ it is smoother and requires more samples than $m = 8q$. This is consistent with Theorem~\ref{thm:recovery_CFC}, where the number of collocation points required to guarantee the sparse recovery depends logarithmically on $d$. These results illustrate that our numerical method lessens the curse of dimensionality in the case for the recovery of sparse solutions.

\subsection{Robustness of the method with respect to the diffusion coefficient $a$}

Proposition~\ref{prop:1-sparse_diffusion} and~\ref{prop:compressible_eta} rely on sufficient conditions on the diffusion coefficient $a$ involving the constants $e_{\bm{0}}$, $\alpha$, $\beta$, and $\gamma$. However, in all the experiments above, we have not discussed the validity of these conditions for the diffusion coefficients $a_1$, $a_2$, and $a_3$. In this section, we numerically investigate the robustness of compressive Fourier collocation with respect to the diffusion coefficient $a$, considering regimes where $a$ satisfies or breaks the sufficient conditions of our theory. 
 
\begin{figure}[t] 
\begin{subfigure}{.5\textwidth}
  \centering
  \includegraphics[width = 0.9\textwidth]{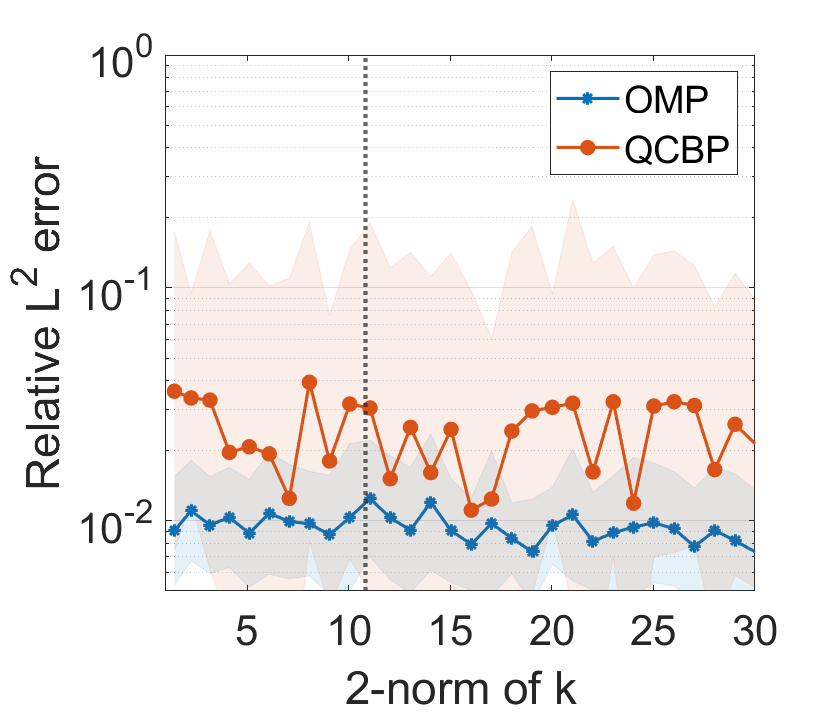}
\end{subfigure}
\begin{subfigure}{.5\textwidth}
  \centering
  \includegraphics[width = 0.9\textwidth]{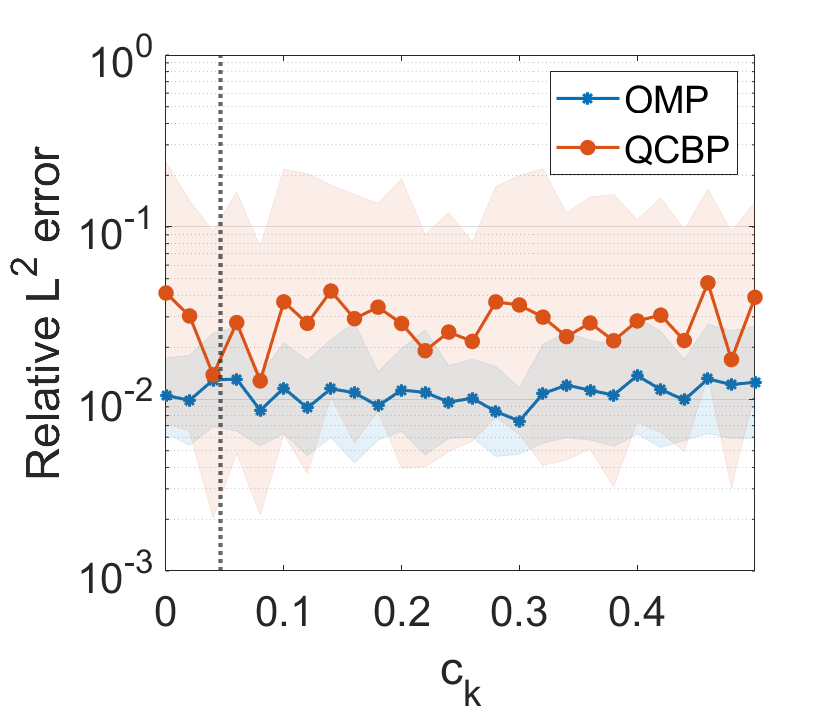}
\end{subfigure}
\caption{(Robustness of the method with respect to the diffusion coefficient; $d=2$, nonsparse solution $u_3$)  Left: (impact of $\bm{k}$ on accuracy) Relative $L^2$-error versus $\|\bm{k}\|_2$.  Right: (impact of $c_{\bm{k}}$ on accuracy) Relative $L^2$-error versus $c_{\bm{k}}$.  } \label{Fig:impact_k_ck}
\end{figure}

We consider a diffusion coefficient of the form \eqref{eq:exa_diffusion_cosine} with dimension $d=2$ and the nonsparse exact solution $u_3$. In this case, recalling Remark~\ref{app:example_diffusion}, condition \eqref{eq:compressible_eta_cond_2} of Proposition~\ref{prop:compressible_eta} is satisfied because $a^* = 0$ and condition \eqref{eq:compressible_eta_cond_1}, involving $e_{\bm{0}}$ and $\beta$, is equivalent to
\begin{equation*}
|c_{\bm{0}}| >  \frac{ \sqrt{1 + (2\pi)^2 \|\bm{k}\|_2^2}}{\sqrt{2}-1} |c_{\bm{k}}|.
\end{equation*} 
To test the robustness of the method, we fix $c_{\bm{0}}$ and study the relative $L^2$-error as a function of $\bm{k}$ and $c_{\bm{k}}$. The results are shown in Fig.~\ref{Fig:impact_k_ck}. Specifically, in the left plot, we test the impact of $\bm{k}$ on the method's accuracy by changing $\bm{k}$ in the diffusion coefficients  of the form $a(\bm{x}) = 1 + 0.02 \cos(2 \pi k_1 x_1) \cos(2 \pi k_2 x_2)$, for $k_1 = 1,\ldots, 30$, and $k_2 = 1$. In the right plot of Fig.~\ref{Fig:impact_k_ck}, we test the impact of $c_{\bm{k}}$ while keeping $\bm{k}$ fixed for diffusion coefficients of the form $a(\bm{x}) = 1 + c_{\bm{k}} \cos{2 \pi x_1} \cos{2 \pi x_2}$, with $\bm{k}=(1,1)$. In both figures, the vertical dotted line shows the value of $\|\bm{k}\|_2$ (or $c_{\bm{k}}$, respectively) for which condition \eqref{eq:compressible_eta_cond_1} becomes an equality. In particular, condition \eqref{eq:compressible_eta_cond_1} is only satisfied on the left of the vertical dotted line. We run 25 random tests for each value of the parameters.

The results of Fig.~\ref{Fig:impact_k_ck} show that the accuracy of the method is not impacted by the validity of condition \eqref{eq:compressible_eta_cond_1}. Hence, it is considerably more robust than what is predicted by the theory.

\section{Conclusions and open problems}
\label{s:conclusions}

We have proposed a new method for the numerical solution of high-dimensional diffusion equations with periodic boundary conditions called compressive Fourier collocation. In Theorem~\ref{thm:recovery_CFC}, we have shown that the proposed method is able to approximate sparse or compressible solutions to the high-dimensional periodic diffusion equation in an accurate and stable way by using a number of collocation points that depends logarithmically on $d$ and that is therefore only mildly affected by the curse of dimensionality. Our numerical experiments in \S\ref{s:numerics} confirm the accuracy and stability of the method up to dimension $d=20$. 

We consider our results very promising and hope they will attract the attention of the scientific machine learning community. However, there are many open problems and potential avenues for future research that remain to be investigated. First, although the compressive Fourier collocation method is able to lessen the curse of dimensionality with respect to the number of collocation points, the number of flops needed in general to recover the solution via OMP or QCBP scales at least linearly in $N=|\Lambda|$, and it is therefore still affected by the curse of dimensionality with respect to the computational complexity. A promising way to address this issue is the combination of the compressive Fourier collocation approach with sublinear sparse recovery techniques, such as those recently proposed in \cite{choi2021sparse,choi2021sparse2}. A significant step forward in this direction has been very recently made by \cite{gross2023sparse}, where the authors propose a new class of a sublinear-time spectral methods for multiscale elliptic PDEs with periodic boundary conditions based on a Fourier-Galerkin discretization and on the sparse Fourier transform. Developing sublinear-time algorithms for compressive Fourier collocation is a key avenue of future work. 

Another important open issue currently under investigation is the extension of the method to PDEs beyond the model problem \eqref{eq:diffusion_eq_periodic_1}--\eqref{eq:diffusion_eq_periodic_2}. A first natural extension is the case of high-dimensional advection-diffusion-reaction equations with periodic boundary conditions, where the PDE operator is of the form $\mathcal{L}[v] = -\nabla \cdot (a \nabla v) + \bm{b} \cdot \nabla v + c v$, with $v\in H^2(\mathbb{T}^d)$ and where $\bm{b} : \mathbb{T}^d \to \mathbb{T}^d$ and $c : \mathbb{T}^d \to \mathbb{C}$. The same type of analysis carried out in this paper should lead to the bounded Riesz property under sufficient conditions on the Fourier coefficients of $\bm{b}$ and $c$, similarly to Proposition~\ref{prop:compressible_eta}. The case of PDEs with nonperiodic boundary conditions is expected to be more challenging. This could be addressed by considering spectral bases based on, e.g., boundary-adapted orthogonal polynomials. However, it is not clear if a rigorous theoretical analysis would be accessible in that case. The ability to deal with nonperiodic boundary conditions would be a key step towards applying compressive spectral collocation to more realistic high-dimensional PDE models.

Another key open issue is the numerical validation of the proposed technique on higher dimensional domains with $d >20$ and its comparison with other collocation-based techniques, such as those based on deep learning mentioned in \S\ref{s:literature}. In this direction, an interesting question would be a comparison in view of the so-called ``curse of high-frequency'', referring to the difficulty of deep neural networks to learn high-frequency information (see \cite{luo2021upper} and references therein). Regarding our theoretical analysis, we observe that conditions \eqref{eq:compressible_eta_cond_1}--\eqref{eq:compressible_eta_cond_2}  on the diffusion coefficient $a$ sufficient to guarantee the bounded Riesz property in Proposition~\ref{prop:compressible_eta} rely on an application of the Gershgorin circle theorem (see proof in \S\ref{app:prop:compressible_eta}). How to improve this proof strategy, likely resulting in weaker conditions on $a$, is an open problem. Moreover, we think that an extremely promising avenue of future work is the combination of our recovery guarantees with recent approximation theory results for deep neural networks (see, e.g., \cite{elbrachter2021}) in order to derive \emph{practical existence theorems} for the numerical solution of high-dimensional PDEs via deep learning, in the same spirit of the results in \cite{adcock2021deep,adcock2021gap} obtained in the context of high-dimensional function approximation. These results provide sufficient conditions on the architecture, the number of data points, and the (regularized) loss function that lead to training deep neural networks with desirable approximation properties. They are based on combining deep learning approximation theory results with recovery theorems for sparse high-dimensional approximation. Hence, our work can be seen as a first important step towards bridging the gap between theory and practice in the emerging area of deep learning-based solution of high-dimensional PDEs.

\section*{Acknowledgments}
The authors gratefully acknowledge the support of the Natural Sciences and Engineering Research Council (NSERC) of Canada through grant RGPIN-2020-06766, the Fonds de Recherche du Québec - Nature et Technologies (FRQNT) through grant number 313276, and the Faculty of Arts and Science of Concordia University. Weiqi Wang is also supported by Concordia University through the Horizon Postdoctoral Fellowship program.

\appendix

\section{Further details and proofs}\label{s:appendix}

\subsection{The Gram matrix}\label{app:Gram}

To show that the system $\{\Phi_{\bm{\nu}}\}_{\bm{\nu} \in \Lambda}$ defined in \eqref{eq:def_Phi} is a Riesz system, it is convenient to consider its \emph{Gram matrix} $G \in \mathbb{C}^{N\times N}$, where $N = |\Lambda|$, defined by
\begin{equation}
\label{eq:def_G}
    G_{\bm{\nu}\bm{\mu}}=\left< \Phi_{\bm{\nu}}, \Phi_{\bm{\mu}}\right>, \quad \forall \bm{\nu},\bm{\mu}\in \Lambda.
\end{equation}
We note in passing that, thanks to the normalization factor $1/\sqrt{m}$ in \eqref{eq:def_A_b}, we have $\mathbb{E}[A^*A] = G$ (this follows from a direct computation and the fact that the random collocation points $\bm{y}_1,\ldots,\bm{y}_m$ are independently and uniformly distributed over $\mathbb{T}^d$). The significance of the Gram matrix $G$ relies on the fact that it yields the following norm equivalence:
\begin{equation}
\label{eq:norm_equivalence}    
\left\| \sum_{\bm{\nu} \in \Lambda} c_{\bm{\nu}} \Phi_{\bm{\nu}} \right\|_{L^2}^2 
= \left< \sum_{\bm{\nu} \in \Lambda} c_{\bm{\nu}} \Phi_{\bm{\nu}},\sum_{\bm{\nu} \in \Lambda} c_{\bm{\nu}} \Phi_{\bm{\nu}}\right>
=\bm{c}^{T}G\bm{c}, \quad \forall \bm{c} \in \mathbb{C}^N.
\end{equation}
Note that $G$ is a Hermitian positive semidefinite matrix. Hence, it has only real nonnegative eigenvalues. The Courant–Fischer–Weyl min-max principle implies that, if $0 < b_{\Phi} \leq B_{\Phi} <\infty$ are such that
\begin{equation}
\label{eq:two-sided_bound}
    b_\Phi \leq \lambda_{\min}(G)\leq \lambda_{\max}(G) \leq B_\Phi,
\end{equation}
then $\{\Phi_{\bm{\nu}}\}_{\bm{\nu}\in\Lambda}$ is a Riesz system with constants $b_{\Phi}$ and $B_{\Phi}$.
Hence, estimating the lower and upper Riesz constants of $\{\Phi_{\bm{\nu}}\}_{\bm{\nu}\in\Lambda}$ corresponds to finding two-sided spectral bounds for the Gram matrix $G$.

Thanks to the Fourier expansion \eqref{eq:eta_expansion} of $a$, it is possible to compute an explicit formula for the entries of the Gram matrix $G$.

\begin{lemma}[Explicit formula for the Gram matrix] \label{lem:expansion_G}
Let $\Lambda \subset \mathbb{Z}^d$ with $\bm{0} \notin \Lambda$ and consider a diffusion coefficient $a \in C^1(\mathbb{T}^d)$ with Fourier expansion as in \eqref{eq:eta_expansion}. Then, the elements of Gram matrix $G$ defined in \eqref{eq:def_G} admit the following explicit formula in terms of the Fourier coefficients $\bm{e} = (e_{\bm{\nu}})_{\bm{\nu} \in \mathbb{Z}^d}$ of $a$:
$$
G_{\bm{\nu} \bm{\mu}} = \sum_{\bm{\tau}\in \mathbb{Z}^d}\sum_{\bm{\tau'}\in \mathbb{Z}^d}  \frac{(\bm{\tau} \cdot \bm{\nu} +  \| \bm{\nu} \|_2^2)}{\| \bm{\nu} \|_2^2} \frac{(\bm{\tau'} \cdot \bm{\mu} + \| \bm{\mu} \|_2^2)}{ \| \bm{\mu} \|_2^2} e_{\bm{\tau}} \bar{e}_{\bm{\tau'}}
\delta_{\bm{\tau}+\bm{\nu},\bm{\tau'}+\bm{\mu}} ,
\quad \forall \bm{\nu}, \bm{\mu} \in \Lambda,
$$
where $\delta_{\bm{\nu},\bm{\mu}}$ denotes the Kronecker delta.
\end{lemma}
\begin{proof} 
Before proving the identity, we note that gradients and Laplacians of the Fourier basis functions defined in \eqref{eq:Fourier_system} can be easily computed as
\begin{equation}
\label{eq:diff_Fourier}
\nabla F_{\bm{\nu}} = (2\pi i \bm{\nu}) F_{\bm{\nu}} \quad \text{and} \quad
\Delta F_{\bm{\nu}} = - 4\pi^2 \|\bm{\nu}\|_2^2 F_{\bm{\nu}}, \quad 
\forall  \bm{\nu} \in \mathbb{Z}^d.
\end{equation}
Moreover, 
\begin{equation}
\label{eq:prod_Fourier}
    F_{\bm{\nu}}F_{\bm{\mu}} = F_{\bm{\nu} +\bm{\mu}},\forall \bm{\nu},\bm{\mu} \in \mathbb{Z}^d.
\end{equation}
Properties \eqref{eq:diff_Fourier} and \eqref{eq:prod_Fourier} are independent of the normalization employed. Hence, they also hold after replacing $F_{\bm{\nu}}$ with the spectral basis functions $\Psi_{\bm{\nu}}$ defined in \eqref{eq:def_Psi}. 

To prove the desired formula for $G_{\bm{\nu}\bm{\mu}}$, we expand the inner product in \eqref{eq:def_G}. Using the above properties, the $L^2$-orthonormality of the Fourier basis $\{F_{\bm{\nu}}\}_{\bm{\nu} \in \mathbb{Z}^d}$, and recalling the expansion \eqref{eq:eta_expansion} of $a$, we see that
\begin{align*}
G_{\bm{\nu} \bm{\mu}} 
& = \langle \nabla a \cdot \nabla \Psi_{\bm{\nu}}
+ a \Delta \Psi_{\bm{\nu}},
 \nabla a \cdot \nabla \Psi_{\bm{\mu}}
+ a \Delta \Psi_{\bm{\mu}}
\rangle \\
& = \sum_{\bm{\tau}}\sum_{\bm{\tau'}}
\langle((2 i \pi \bm{\tau}) \cdot (2 i \pi \bm{\nu}) - 4 \pi^2
\| \bm{\nu} \|_2^2)  e_{\bm{\tau}} F_{\bm{\tau}}\Psi_{\bm{\nu}},  ((2 i \pi \bm{\bm{\tau}}') \cdot (2 i \pi \bm{\mu}) - 4 \pi^2 \| \bm{\mu} \|_2^2) 
\bar{e}_{\bm{\tau'}} F_{\bm{\tau'}}\Psi_{\bm{\mu}}\rangle \\
& = \sum_{\bm{\tau}}\sum_{\bm{\tau'}}  \frac{(4 \pi^2\bm{\tau} \cdot \bm{\nu} + 4 \pi^2 \| \bm{\nu} \|_2^2)}{4 \pi^2 \| \bm{\nu} \|_2^2} \frac{(4 \pi^2\bm{\tau'} \cdot \bm{\mu} + 4 \pi^2 \| \bm{\mu} \|_2^2)}{4 \pi^2 \| \bm{\mu} \|_2^2} e_{\bm{\tau}} \bar{e}_{\bm{\tau'}}
\langle  F_{\bm{\tau}} F_{\bm{\nu}},  F_{\bm{\tau'}} F_{\bm{\mu}}\rangle \\
& =  \sum_{\bm{\tau}}\sum_{\bm{\tau'}}  \frac{(\bm{\tau} \cdot \bm{\nu} +  \| \bm{\nu} \|_2^2)}{\| \bm{\nu} \|_2^2} \frac{(\bm{\tau'} \cdot \bm{\mu} + \| \bm{\mu} \|_2^2)}{ \| \bm{\mu} \|_2^2} e_{\bm{\tau}} \bar{e}_{\bm{\tau'}}
\langle  F_{\bm{\tau}+\bm{\nu}},  F_{\bm{\tau'}+\bm{\mu}}\rangle \\
& = \sum_{\bm{\tau}}\sum_{\bm{\tau'}}  \frac{(\bm{\tau} \cdot \bm{\nu} +  \| \bm{\nu} \|_2^2)}{\| \bm{\nu} \|_2^2} \frac{(\bm{\tau'} \cdot \bm{\mu} + \| \bm{\mu} \|_2^2)}{ \| \bm{\mu} \|_2^2} e_{\bm{\tau}} \bar{e}_{\bm{\tau'}}
\delta_{\bm{\tau}+\bm{\nu},\bm{\tau'}+\bm{\mu}}, 
\end{align*}
where all the summations are over $\bm{\tau},\bm{\tau}' \in \mathbb{Z}^d$.
\end{proof}

\subsection{Cardinality bound for the hyperbolic cross} \label{app:HCbounds}

We illustrate how to obtain the cardinality bound in \eqref{eq:cardinality_bound_Lambda}. Applying \citet[Theorem 3.7]{chernov2016new} with $s = d$, $a = 1$, and $T = n$ yields the bound
$$
|\Lambda^{\mathsf{HC}}_{d,n-1}| < 2 \delta^{-1} n^{1 + 2/\delta} (1-\delta)^{-2d/\delta},
$$
for any $0 < \delta < 1$. In particular, choosing $\delta = 1/2$ leads to 
$$
|\Lambda^{\mathsf{HC}}_{d,n-1}| < 4 n^5 16^d.
$$
A second bound can be found in the proof of \citet[Theorem 4.9]{kuhn2015approximation}. Therein, the cardinality of the hyperbolic cross is denoted by $C(r,d)$. Setting $r = n$ leads to the bound
$$
|\Lambda^{\mathsf{HC}}_{d,n-1}| \leq \mathrm{e}^2 n^{2+\log_2(d)}.
$$
Combining the above bounds yields \eqref{eq:cardinality_bound_Lambda}.

\subsection{The Gershgorin circle theorem}
\label{app:Gershgorin}

We recall the Gershgorin circle theorem, whose proof can be found, e.g., in \citet[Theorem 6.1.1]{horn2012matrix}.
\begin{theorem}[Gershgorin circle theorem] 
Let $M \in \mathbb{C}^{N \times N}$. Define the Gershgorin disc of the $i$-th row $D(M_{ii},R_i)$ as the closed disc centered at $M_{ii}$ with radius $R_i = \sum_{j \neq i}|M_{ij}|$. Then every eigenvalue of $M$ lies in at least one of the Gershgorin discs $D(M_{ii},R_i)$.
\end{theorem}
In our proofs, we will repeatedly use the following immediate consequence of the Gershgorin circle theorem. 
\begin{corollary}[Gershgorin circle theorem for Hermitian matrices] 
\label{Gershgorin}
Let $M \in \mathbb{C}^{N \times N}$ a Hermitian matrix. Then, all eigenvalues of $M$ lie in the real interval
\begin{align*}
\left[\min_{i\in[N]}\left\{M_{ii} -  \sum_{j \neq i}|M_{ij}| \right\},
%\leq \lambda \leq 
\max_{i\in[N]}\left\{M_{ii} + \sum_{j \neq i}|M_{ij}| \right\}\right].
\end{align*}
\end{corollary}

\subsection{Proof of Proposition~\ref{prop:1-sparse_diffusion}}
\label{app:prop:1-sparse_diffusion}

The proof has two main steps. First, we show the Riesz property. Then, we prove that the Riesz system is bounded.

\paragraph{Step 1: Riesz property.} As explained in Appendix~\ref{app:Gram}, it is sufficient to find a two-sided spectral bound for the Gram matrix $G$ of the form \eqref{eq:two-sided_bound}. By assumption, $\bm{e} = (e_{\bm{\nu}})_{\bm{\nu}\in\mathbb{Z}^d}$ has only two non-zero terms: $e_{\bm{0}}$ and $e_{\bm{\nu*}}$. Then, using Lemma~\ref{lem:expansion_G}, the only nonzero entries in the $\bm{\nu}$-th row of $G$ are
\begin{align*}
    G_{\bm{\nu} \bm{\nu}} & = \left|e_{\bm{0}}\right|^2+\left( \frac{\bm{\nu}^* \cdot \bm{\nu} +  \| \bm{\nu} \|_2^2}{\| \bm{\nu} \|_2^2} \right)^2 \left|e_{\bm{\nu}^*}\right|^2 =
    \left|e_{\bm{0}}\right|^2+\left( 1+ \frac{\bm{\nu}^* \cdot \bm{\nu}}{\bm{\nu} \cdot \bm{\nu}}\right)^2 |e_{\bm{\nu}^*}|^2.\\
    G_{\bm{\nu},\bm{\nu}-\bm{\nu}^*} & = \frac{\bm{\nu}^* \cdot (\bm{\nu}-\bm{\nu}^*) + \| \bm{\nu}-\bm{\nu}^* \|_2^2}{ \| \bm{\nu}-\bm{\nu}^* \|_2^2} e_{\bm{0}} \bar{e}_{\bm{\nu}^*}
    = \left( \frac{\bm{\nu}^* \cdot (\bm{\nu}-\bm{\nu}^* )}{ \| \bm{\nu}-\bm{\nu}^* \|_2^2} + 1 \right) e_{\bm{0}} \bar{e}_{\bm{\nu}^*}. \\
    G_{\bm{\nu},\bm{\nu}+\bm{\nu}^*} & = \frac{\bm{\nu}^* \cdot \bm{\nu} + \| \bm{\nu} \|_2^2}{ \| \bm{\nu} \|_2^2} \bar{e}_{\bm{0}} e_{\bm{\nu}^*} 
    = \left( \frac{\bm{\nu}^* \cdot \bm{\nu}}{\| \bm{\nu} \|_2^2}+1\right) \bar{e}_{\bm{0}} e_{\bm{\nu}^*}.
\end{align*}
Recalling that $\alpha = |e_{\bm{\nu}^*}| (2\|\bm{\nu}^*\|_2+3)$, and using the Cauchy-Schwarz inequality combined with the fact that $\|\bm{\nu}^*\|_2 \geq 1$,  we derive the following two-sided bound for the diagonal entries of $G$:
$$
\left|e_{\bm{0}}\right|^2 
\leq G_{\bm{\nu} \bm{\nu}}
\leq \left|e_{\bm{0}}\right|^2+ \left( \frac{1+ \|\bm{\nu}^*\|_2}{2\|\bm{\nu}^*\|_2 +3 }\right)^2 \alpha^2 
\leq \left|e_{\bm{0}}\right|^2 + \frac{\alpha^2}{4}.
$$
Moreover, the sum of the absolute values of off-diagonal entries of $G$ can be bounded as follows:
\begin{align*}
    \left|G_{\bm{\nu},\bm{\nu}-\bm{\nu}^*}\right| + \left|G_{\bm{\nu}, \bm{\nu}+\bm{\nu}^*}\right| 
    & \leq  \left(\frac{\| \bm{\nu}^* \|_2}{ \| \bm{\nu}-\bm{\nu}^* \|_2} + \frac{\| \bm{\nu}^* + \bm{\nu} \|_2 }{\| \bm{\nu} \|_2} +2 \right) \cdot \frac{\left|e_{\bm{0}}\right|\alpha}{2\|\bm{\nu}^*\|_2+3} \\
    & \leq \left( \| \bm{\nu}^* \|_2 + \frac{\| \bm{\nu}^*\|_2}{\| \bm{\nu} \|_2}  + 3  \right)\cdot \frac{\left|e_{\bm{0}}\right|\alpha}{2\|\bm{\nu}^*\|_2+3} \\
    & \leq (2\| \bm{\nu}^* \|_2+3) \cdot \frac{\left|e_{\bm{0}}\right|\alpha}{2\|\bm{\nu}^*\|_2+3}
     = \left|e_{\bm{0}}\right| \alpha, 
\end{align*}
where we used the Cauchy-Schwarz inequality, the definition of $\alpha$, the triangle inequality, combined with the facts that $\|\bm{\nu}^*\|_2 \geq 1$ (by assumption) and that $\|\bm{\nu}\|_2\geq1$ and $ \|\bm{\nu}-\bm{\nu}^*\|_2 \geq 1$ since $\bm{0} \notin \Lambda$. 

Finally, applying Corollary \ref{Gershgorin} we have 
\begin{align*}
    \lambda_{\min}(G)
    & \geq G_{\bm{\nu} \bm{\nu}} - \left|G_{\bm{\nu},\bm{\nu}-\bm{\nu}^*}\right| - \left|G_{\bm{\nu}, \bm{\nu}+\bm{\nu}^*}\right| 
    \geq \left|e_{\bm{0}}\right|^2- \left|e_{\bm{0}}\right| \alpha =: b_{\Phi},\\
    \lambda_{\min}(G)
    & \leq G_{\bm{\nu} \bm{\nu}} + \left|G_{\bm{\nu},\bm{\nu}-\bm{\nu}^*}\right| + \left|G_{\bm{\nu}, \bm{\nu}+\bm{\nu}^*}\right| 
    \leq \left|e_{\bm{0}}\right|^2 + 
    \left|e_{\bm{0}}\right| \alpha + \frac{\alpha^2}{4}  = \left(|e_{\bm{0}}| + \frac{\alpha}{2}\right)^2 =: B_{\Phi},
\end{align*}
which is the desired two-sided bound.

\paragraph{Step 2: Boundedness.} To determine the uniform bound for the system, recalling \eqref{eq:diff_Fourier} we compute
\begin{align*}
    \Phi_{\bm{\nu}} & =
    \frac{1}{4\pi^2 \| \bm{\nu}\|_2^2}\left(-\nabla\cdot((e_{\bm{0}}+e_{\bm{\nu}^*}F_{\bm{\nu}^*})(2\pi i \bm{\nu})F_{\bm{\nu}}) \right)\\
    & = \frac{1}{4\pi^2 \| \bm{\nu}\|_2^2} \left( 
    4\pi^2 \| \bm{\nu}\|_2^2 e_{\bm{0}} F_{\bm{\nu}} 
    -e_{\bm{\nu}}^*(2\pi i \bm{\nu})  \cdot \nabla
    F_{\bm{\nu} + \bm{\nu}^*}
    \right) \\
    & = e_{\bm{0}} F_{\bm{\nu}} +
    e_{\bm{\nu}^*} F_{\bm{\nu}+\bm{\nu}^*}\frac{\bm{\nu}\cdot(\bm{\nu}+\bm{\nu}^*)}{\| \bm{\nu}\|_2^2}.
    %& = e_\bm{0} F_{\bm{\nu}} +
    %e_\bm{\nu^*} F_{\bm{\nu}+\bm{\nu^*}} + e_\bm{\nu^*} F_{\bm{\nu}+\bm{\nu^*}} \frac{\bm{\nu} \cdot \bm{\nu^*}}{\| \bm{\nu}\|_2^2}
\end{align*}
As a consequence, we obtain the the bound
$$
\left\| \Phi_{\bm{\nu}} \right\|_{L^\infty} 
\leq \left| e_{\bm{0}} \right| +
\left| e_{\bm{\nu}^*} \right|\left(1+\frac{|\bm{\nu} \cdot \bm{\nu}^*|}{\| \bm{\nu}\|_2^2}\right).
$$
Finally, using Cauchy-Schwarz, the definition of $\alpha$, and the fact that $\|\bm{\nu}\|_2 \geq 1$, we see that
\begin{align*}
\left\| \Phi_{\bm{\nu}} \right\|_{L^\infty}  
 \leq \left| e_{\bm{0}} \right| +
    \left| e_{\bm{\nu}^*} \right|\left(1+\frac{\|\bm{\nu}^*\|_2}{\| \bm{\nu}\|_2}\right) 
     \leq \left| e_{\bm{0}} \right| +
    \left( \frac{\alpha}{2\|\bm{\nu}^*\|_2+3} \right)\left(1+\|\bm{\nu}^*\|_2\right) 
     \leq \left| e_{\bm{0}} \right| + \frac{\alpha}{2}. 
\end{align*}
This concludes the proof.

\subsection{Proof of Proposition~\ref{prop:compressible_eta}}
\label{app:prop:compressible_eta}

The proof is structured in two main parts. First, in \S\ref{s:eta*=0} we consider the simpler case where $a^* = 0$ (and, hence, $a = a_t$). Then,  in \S\ref{s:eta*neq0} we extend the proof to the general case $a^* \neq 0$. 

\subsubsection{The case $a^* = 0$}
\label{s:eta*=0}

\paragraph{Case $a^* = 0$, Step 1: Riesz property.} Similarly to the case of Proposition~\ref{prop:1-sparse_diffusion}, we find lower and upper Riesz constants $b_{\Phi}$ and $B_{\Phi}$ by establishing a two-sided spectral bound for the Gram matrix $G$. Using Lemma~\ref{lem:expansion_G}, the diagonal entries of $G$ are given by
$$
G_{\bm{\nu} \bm{\nu}}  
=  \sum_{\bm{\tau} \in T \cup \{\bm{0}\}}  \frac{(\bm{\tau} \cdot \bm{\nu} +  \| \bm{\nu} \|_2^2)^2}{\| \bm{\nu} \|_2^4}  \lvert e_{\bm{\tau}} \rvert^2, \quad \forall \bm{\nu} \in \Lambda. 
$$
Using the Cauchy-Schwarz inequality and the fact that $\|\bm{\nu}\|_2\geq 1$, we see that %$ G_{\bm{\nu} \bm{\nu}}  $ has the upper bound and lower bound using $\lvert \bm{\tau} \cdot \bm{\nu}\rvert \leq \| \bm{\tau} \|_2 \cdot \| \bm{\nu} \|_2$
\begin{equation} 
\label{eq:G_diag}
\lvert e_{\bm{0}} \rvert^2 
\leq  G_{\bm{\nu} \bm{\nu}}   \leq  \sum_{\bm{\tau} \in T \cup \{\bm{0}\}} \left(\frac{\| \bm{\tau} \|_2}{\| \bm{\nu} \|_2} + 1 \right)^2 \lvert e_{\bm{\tau}} \rvert^2
\leq \lvert e_{\bm{0}} \rvert^2+ \sum_{\bm{\tau} \in T } \left(\| \bm{\tau} \|_2 + 1 \right)^2 \lvert e_{\bm{\tau}} \rvert^2 
    \leq \left\| a_t \right\|^2_{H^1}.
\end{equation}
The above inequality can be proved as follows. Using the definition of $H^1$-norm, the differentiation properties \eqref{eq:diff_Fourier}, and the fact that the Fourier system $\{F_{\bm{\nu}}\}_{\bm{\nu} \in \mathbb{Z}^d}$ is $L^2$-orthonormal, we obtain 
\begin{align*}
     \left\| a_t \right\|^2_{H^1} 
     & =   \left\| a_t \right\|^2_{L^2}
 + \sum^{d}_{l=1}  \left\| \frac{\partial a_t}{\partial x_l} \right\|^2_{L^2}\\
     &= \sum_{\bm{\tau}\in T \cup \{\bm{0}\}}|e_{\bm{\tau}}|^2 + \sum^{d}_{l=1}\sum_{\bm{\tau}\in T \cup \{\bm{0}\}}|2\pi e_{\bm{\tau}} \bm{\tau}_l|^2\\
     %& = \sum_{\bm{\tau}\in T \cup \{\bm{0}\}}(1+(2\pi)^2 \left\|\bm{\tau}\right\|_2^2)|e_\bm{\tau}|^2\\
     & = |e_{\bm{0}}|^2 + \sum_{\bm{\tau} \in T}\left(1+(2\pi)^2 \left\|\bm{\tau}\right\|_2^2\right)|e_{\bm{\tau}}|^2\\
     & \geq  \lvert e_{\bm{0}} \rvert^2+ \sum_{\bm{\tau} \in T} \left(\| \bm{\tau} \|_2 + 1 \right)^2 \lvert e_{\bm{\tau}} \rvert^2, 
\end{align*}
which proves \eqref{eq:G_diag}.

To apply Gershgorin circle theorem, we now bound the sum of all off-diagonal entries in the $\bm{\nu}$-th row of $G$. Using Lemma~\ref{lem:expansion_G} again, the definition of the Kronecker delta, the Cauchy-Schwarz inequality, and the fact that $\bm{0} \notin \Lambda$ yields
$$
\sum_{\bm{\mu} \in \Lambda \setminus\{ \bm{\nu}\}} \lvert G_{\bm{\nu} \bm{\mu}} \rvert
\leq \sum_{\bm{\tau} \in T \cup \{\bm{0}\}} \left(\| \bm{\tau} \|_2 + 1 \right) \lvert e_{\bm{\tau}} \rvert \sum_{\bm{\mu} \in \Lambda \setminus\{ \bm{\nu}\}} \left(\| \bm{\tau} +\bm{\nu} -\bm{\mu}\|_2 + 1 \right) \lvert e_{\bm{\bm{\tau} +\bm{\nu} -\bm{\mu}}} \rvert.  
$$
Substituting $\bm{\tau}'=\bm{\tau}+\bm{\nu}-\bm{\mu}$ (which implies $\bm{\tau}' \neq \bm{\tau}$), recalling that $\bm{e}$ is supported on $T$, and separating the $e_{\bm{0}}$ term, we obtain
\begin{align*}
    \sum_{\bm{\mu} \in \Lambda \setminus\{ \bm{\nu}\}} \lvert G_{\bm{\nu} \bm{\mu}} \rvert
    & \leq \sum_{\bm{\tau} \in T \cup \{\bm{0}\}} \left(\| \bm{\tau} \|_2 + 1 \right) \lvert e_{\bm{\tau}} \rvert \sum_{\bm{\tau'}  \in T \cup \{\bm{0}\} \setminus\{ \bm{\tau}\}} \left(\| \bm{\tau'} \|_2 + 1 \right) \lvert e_{\bm{\tau'}} \rvert \\
    & = \lvert e_{\bm{0}} \rvert \sum_{\bm{\tau'}  \in T } \left(\| \bm{\tau'} \|_2 + 1 \right) \lvert e_{\bm{\tau'}} \rvert + \sum_{\bm{\tau}  \in T } \left(\| \bm{\tau} \|_2 + 1 \right) \lvert e_{\bm{\tau}} \rvert \sum_{\bm{\tau'}  \in T \setminus\{ \bm{\tau}\}} \left(\| \bm{\tau'} \|_2 + 1 \right) \lvert e_{\bm{\tau'}} \rvert \\
    & \leq \lvert e_{\bm{0}} \rvert \sum_{\bm{\tau'}  \in T } \left(\| \bm{\tau'} \|_2 + 1 \right) \lvert e_{\bm{\tau'}} \rvert + \sum_{\bm{\tau}  \in T} \left(\| \bm{\tau} \|_2 + 1 \right) \lvert e_{\bm{\tau}} \rvert \sum_{\bm{\tau'} \in T} \left(\| \bm{\tau'} \|_2 + 1 \right) \lvert e_{\bm{\tau'}} \rvert \\
    & = 2 \lvert e_{\bm{0}} \rvert \sum_{\bm{\tau'}  \in T } \left(\| \bm{\tau'} \|_2 + 1 \right) \lvert e_{\bm{\tau'}} \rvert + \left( \sum_{\bm{\tau}  \in T } \left(\| \bm{\tau} \|_2 + 1 \right) \lvert e_{\bm{\tau}} \rvert \right)^2.
\end{align*}
%Let $p=\sum_{\bm{\tau'}\neq 0} \left(\| \bm{\tau'} \| + 1 \right) \lvert e_{\bm{\tau'}} \rvert$. 
Applying the Cauchy–Schwarz inequality and using that $|T| \leq t$, we obtain the bound 
\begin{equation}
\label{eq:Bound_beta}
    \sum_{\bm{\tau} \in T } \left(\| \bm{\tau} \|_2 + 1 \right) \lvert e_{\bm{\tau}} \rvert 
    \leq \sqrt{t}
    \left( \sum_{\bm{\tau}\in T} \left(\| \bm{\tau} \|_2 + 1 \right)^2 \lvert e_{\bm{\tau}} \rvert ^2 \right)^{1/2}
     \leq \sqrt{t} \left( \left\| a \right\|^2_{H^1} - |e_{\bm{0}}|^2 \right)^{1/2} = \beta.
\end{equation}
Combining the above inequalities yields
\begin{equation} \label{eq:G_off_diag}
    \sum_{\bm{\mu} \in \Lambda \setminus \{\bm{\nu}\} } \lvert G_{\bm{\nu} \bm{\mu}} \rvert %\leq 2\lvert e_{\bm{0}} \rvert p+p^2 
    \leq 2|e_{\bm{0}}| \beta+\beta^2.
\end{equation}

Finally, applying the Gershgorin circle theorem on $G$ combining \eqref{eq:G_diag} and \eqref{eq:G_off_diag}, we obtain the Riesz constants 
\begin{equation} \label{eq:eta_sparse_conclusion}
b_{\Phi} = |e_{\bm{0}}|^2- 2|e_{\bm{0}}| \beta - \beta^2 \quad \text{and} \quad 
B_{\Phi} = \left\|a_t \right\|_{H^1}^2 + 2|e_{\bm{0}}| \beta+ \beta^2.
\end{equation}

\paragraph{Case $a^* = 0$, Step 2: Boundedness.} To determine the uniform bound for the system, we use again properties \eqref{eq:diff_Fourier} and \eqref{eq:prod_Fourier} and compute
\begin{align*}
    \Phi_{\bm{\nu}} 
    %& =     -\nabla\cdot\left(a \nabla \Psi_{\bm{\nu}} \right)
    &= - a_t \Delta \Psi_{\bm{\nu}}- \nabla a_t \cdot\nabla \Psi_{\bm{\nu}} \\ 
    & =  a_t F_{\bm{\nu}} - \left( \sum_{\bm{\tau} \in T \cup \{\bm{0}\}} 2\pi i \bm{\tau}
    e_{\bm{\tau}} F_{\bm{\tau}}\right) \cdot
    \left( \frac{1}{4\pi^2 \| \bm{\nu}\|_2^2} 2\pi i \bm{\nu} F_{\bm{\nu}}\right) \\
    & = \sum_{\bm{\tau}\in T \cup \{\bm{0}\}} e_{\bm{\tau}} F_{\bm{\tau}+\bm{\nu}} +
    \sum_{\bm{\tau}\in T \cup \{\bm{0}\}} \frac{\bm{\tau} \cdot \bm{\nu}}{\| \bm{\nu}\|_2^2} e_{\bm{\tau}} F_{\bm{\tau}+\bm{\nu}} \\
    &= \sum_{\bm{\tau}\in T \cup \{\bm{0}\}}\left(1+\frac{\bm{\tau} \cdot \bm{\nu}}{\| \bm{\nu}\|_2^2}\right) e_{\bm{\tau}} F_{\bm{\tau}+\bm{\nu}}.
\end{align*}
Using that $|\bm{\tau} \cdot \bm{\nu}| \leq \| \bm{\tau} \|_2 \| \bm{\nu} \|_2 $, $\| \bm{\nu} \|_2 \geq 1 $, $\|F_{\bm{\nu}}\|_{L^{\infty}}=1$, and \eqref{eq:Bound_beta}, we obtain 
\begin{equation} \label{eq:eta*=0_infinite_norm}
\left\| \Phi_{\bm{\nu}}
\right\|_{L^\infty} 
\leq\sum_{\bm{\tau}\in T \cup \{\bm{0}\}}\left(1+\frac{|\bm{\tau} \cdot \bm{\nu}|}{\| \bm{\nu}\|_2^2}\right) |e_{\bm{\tau}}|
 \leq |e_{\bm{0}}| + \sum_{\bm{\tau}\in T} \left(1+ \| \bm{\tau} \|_2^2 \right) \lvert e_{\bm{\tau}} \rvert 
\leq |e_{\bm{0}}| + \beta.
\end{equation}
This concludes the proof for $a^* =0$

\subsubsection{The case $a^* \neq 0$}
\label{s:eta*neq0}

We now consider the general case $a^* \neq 0$, where $a = a_t + a^*$.

\paragraph{Case $a^*\neq 0$, Step 1: Riesz property.} Let $\bm{z} \in \ell^2(\Lambda, \mathbb{C}) \cong \mathbb{C}^N$. Keeping the Definition~\ref{def:Riesz} of Bounded Riesz System in mind and using the triangle inequality of the $L^2$-norm, we estimate
\begin{align}
\nonumber
\left\| \sum_{\bm{\nu} \in \Lambda} z_{\bm{\nu}}
    \Phi_{\bm{\nu}} \right\|_{L^2}
    &=
    \left\| \sum_{\bm{\nu} \in \Lambda} z_{\bm{\nu}}
    \nabla \cdot (a \nabla\Psi_{\bm{\nu}}) \right\|_{L^2}\\
    &\geq \left\| \sum_{\bm{\nu} \in \Lambda} z_{\bm{\nu}}
    \nabla \cdot (a_t \nabla \Psi_{\bm{\nu}}) \right\|_{L^2} -
    \left\| \sum_{\bm{\nu} \in \Lambda} z_{\bm{\nu}}
    \nabla \cdot (a^* \nabla \Psi_{\bm{\nu}}) \right\|_{L^2}
    \label{eq:eta*_L2_lower_bound}
\end{align}
and, similarly,
\begin{equation}
\label{eq:eta*_L2_upper_bound}
   \left\| \sum_{\bm{\nu} \in \Lambda} z_{\bm{\nu}}
    \Phi_{\bm{\nu}} \right\|_{L^2}
    %= \left\| \sum_{\bm{\nu} \in \Lambda} z_{\bm{\nu}} \nabla \cdot (a \nabla\Psi_{\bm{\nu}}) \right\|_{L^2}
    \leq \left\| \sum_{\bm{\nu} \in \Lambda} z_{\bm{\nu}}
    \nabla \cdot (a_t \nabla \Psi_{\bm{\nu}}) \right\|_{L^2} +
    \left\| \sum_{\bm{\nu} \in \Lambda} z_{\bm{\nu}}
    \nabla \cdot (a^* \nabla \Psi_{\bm{\nu}}) \right\|_{L^2}.
\end{equation}
We computed lower and upper Riesz constants in the case of a $(t+1)$-sparse diffusion coefficient $a = a_t$ in \S\ref{s:eta*=0}. This leads to two-sided bounds for the first terms in the right-hand sides of the above inequalities. Therefore, we discuss the term involving the term $a^*$. Applying the triangle inequality again, we see that
\begin{align}
\nonumber
    \left\| \sum_{\bm{\nu} \in \Lambda} z_{\bm{\nu}}
    \nabla \cdot (a^* \nabla \Psi_{\nu}) \right\|_{L^2} & = 
    \left\| \sum_{\bm{\nu} \in \Lambda} z_{\bm{\nu}}
    \left( \nabla a^* \cdot \nabla \Psi_{\bm{\nu}} + a^* \Delta \Psi_{\bm{\nu}} \right) \right\|_{L^2} \\
    & \leq \left\| \sum_{\bm{\nu} \in \Lambda} z_{\bm{\nu}}
    \nabla a^* \cdot \nabla \Psi_{\bm{\nu}} \right\|_{L^2}
    + \left\| \sum_{\bm{\nu} \in \Lambda} z_{\bm{\nu}}
    a^* \Delta \Psi_{\bm{\nu}} \right\|_{L^2}.
    \label{eq:eta*_triangle_inequality}
\end{align}
For the first term of \eqref{eq:eta*_triangle_inequality}, we use properties \eqref{eq:diff_Fourier} to obtain
$$
\left| \nabla a^* \cdot \nabla \Psi_{\bm{\nu}} \right|  = 
\left| \sum_{l=1}^d \left(\frac{\partial a^*}{\partial x_l}\right)
\left(\frac{\partial \Psi_{\bm{\nu}}}{\partial x_l}\right) \right| 
= \left| \sum_{l=1}^d \left(\frac{\partial a^*}{\partial x_l} 
\frac{2 \pi i \bm{\nu}_i}{4\pi^2 \| \bm{\nu} \|_2^2}\right) F_{\bm{\nu}} \right|.
$$
Then, taking the absolute value of the Fourier coefficients and applying Cauchy–Schwarz yields
$$
    \left| \nabla a^* \cdot \nabla \Psi_{\bm{\nu}} \right|
     \leq \sum_{l=1}^d \left| \frac{\partial a^*}{\partial x_l}  \right| 
    \left| \frac{2 \pi i \bm{\nu}_l}{4\pi^2 \| \bm{\nu} \|_2^2}  \right| 
     \leq  \|\nabla a^*\|_2 
    \sqrt{\sum_{l=1}^d\left| \frac{2 \pi i \bm{\nu}_l}{4\pi^2 \| \bm{\nu} \|_2^2}  \right|^2} 
     = \frac{1}{2\pi \| \bm{\nu} \|_2} \|\nabla a^*\|_2,
$$
where $\|\nabla a^*\|_2$ denotes the function defined by $\bm{x} \mapsto \|\nabla a^*(\bm{x})\|_2$, for every $\bm{x} \in \mathbb{T}^d$. Using the above inequality, we estimate the $L_2$-norm of the first term in \eqref{eq:eta*_triangle_inequality} as
\begin{equation}
\label{eq:eta*_bound_term_1}
\left\| \sum_{\bm{\nu} \in \Lambda} z_{\bm{\nu}}
\nabla a^* \cdot \nabla \Psi_{\bm{\nu}} \right\|_{L^2}  \leq
\left\| \sum_{\bm{\nu} \in \Lambda} \left| z_{\bm{\nu}} \right|
\left| \nabla a^* \cdot \nabla \Psi_{\bm{\nu}} \right| \right\|_{L^2} 
\leq \left\| \|\nabla a^*\|_2 \sum_{\bm{\nu} \in \Lambda} \frac{1}{2\pi \| \bm{\nu} \|_2}
\left| z_{\bm{\nu}} \right| \right\|_{L^2}.
\end{equation}
%Notice that $| z_{\bm{\nu}} |/2\pi \| \bm{\nu} \|_2$ is a function of $\bm{\nu}$ and $\|\nabla a^*\|_2 = |a^*|_{H^1}$ is independent to $\bm{\nu}$.
Using the definition of $H^1$-norm and the inequality $\|\bm{z}\|_1 \leq \sqrt{N} \|\bm{z}\|_2$, we see that
\begin{align}
\nonumber
    \left\| \|\nabla a^*\|_2  \sum_{\bm{\nu} \in \Lambda} \frac{1}{2\pi \| \bm{\nu} \|_2}
    \left| z_{\bm{\nu}} \right| \right\|_{L^2}
     & = \left(\sum_{\bm{\nu} \in \Lambda} \frac{1}{2\pi \| \bm{\nu} \|_2}    \left| z_{\bm{\nu}} \right| \right)   
    |a^*|_{H^1}\\
\nonumber
    &\leq \frac{1}{2\pi} \left(\sum_{\bm{\nu} \in \Lambda} 
    \left| z_{\bm{\nu}} \right| \right)  
    |a^*|_{H^1} \\
\nonumber
    & = \frac{1}{2\pi} \left\| \bm{z} \right\| _1 
    |a^*|_{H^1}\\
\label{eq:eta*_bound_term_1b}
    & \leq \frac{\sqrt{N}}{2\pi} \left\| \bm{z} \right\| _2 |a^*|_{H^1}. 
\end{align}
Now, using the Cauchy–Schwarz inequality, properties \eqref{eq:diff_Fourier} and the $L^2$-orthonormality of the Fourier system, the second term in \eqref{eq:eta*_triangle_inequality} can be bounded as
\begin{equation}
\label{eq:eta*_bound_term_2}
    \left\| \sum_{\bm{\nu} \in \Lambda} z_{\bm{\nu}}
    a^* \Delta \Psi_{\bm{\nu}} \right\|_{L^2} =
    \left\| \sum_{\bm{\nu} \in \Lambda} z_{\bm{\nu}}
    a^* F_{\bm{\nu}} \right\|_{L^2}
     \leq \left\| a^* \right\|_{L^2} 
    \left\| \sum_{\bm{\nu} \in \Lambda} z_{\bm{\nu}} F_{\bm{\nu}} \right\|_{L^2} 
    \leq \left\| a^* \right\|_{L^2} \left\| \bm{z} \right\| _2.
\end{equation}
Combining \eqref{eq:eta*_triangle_inequality}, \eqref{eq:eta*_bound_term_1}, \eqref{eq:eta*_bound_term_1b}, and \eqref{eq:eta*_bound_term_2} yields
\begin{equation}
    \left\| \sum_{\bm{\nu} \in \Lambda} z_{\bm{\nu}}
    \nabla \cdot (a^* \nabla \Psi_{\bm{\nu}}) \right\|_{L^2} \leq
    \left(\frac{\sqrt{N}}{2\pi}|a^*|_{H^1} + \left\| a^* \right\|_{L^2} \right)  
    \left\| \bm{z} \right\| _2 
    = \gamma \left\| \bm{z} \right\| _2 \label{eq:eta*_bound}
\end{equation}
Now, note that the conclusion \eqref{eq:eta_sparse_conclusion} of \S\ref{s:eta*=0} (Step 1) implies 
\begin{equation*}
    \left( |e_{\bm{0}}|^2- 2|e_{\bm{0}}| \beta - \beta^2 \right) \|\bm{z}\|^2_2 \leq \left\| \sum_{\bm{\nu} \in \Lambda}  z_{\bm{\nu}}
    \nabla (a_t \nabla \Psi_{\bm{\nu}}) \right\|^2_{L^2} \leq \left( \left\|a_t \right\|_{H^1}^2 + 2 |e_{\bm{0}}| \beta+ \beta^2 \right) \|\bm{z}\|^2_2.   
\end{equation*}
Combining the above two-sided bound with \eqref{eq:eta*_L2_lower_bound}, \eqref{eq:eta*_L2_upper_bound}, and \eqref{eq:eta*_bound} we finally obtain the desired estimates
\begin{equation*}
b_{\Phi} = \left( \sqrt{ |e_{\bm{0}}|^2- 2 |e_{\bm{0}}| \beta - \beta^2 } -\gamma \right)^2 \quad\text{and}\quad 
B_{\Phi} = \left( \sqrt{ \left\|a_t \right\|_{H^1}^2 + 2|e_{\bm{0}}| \beta+ \beta^2} +\gamma \right)^2.
\end{equation*}

\paragraph{Case $a^*\neq 0$, Step 2: Boundedness.} To determine the $L^\infty$-norm bound for the system, using properties \eqref{eq:diff_Fourier} and the fact that $\bm{0} \notin \Lambda$, we estimate, for any $\bm{\nu} \in \Lambda$
\begin{align*}
     \left|
    \Phi_{\bm{\nu}} \right|
    & = \left| - a^* \Delta \Psi_{\bm{\nu}}- \nabla a^* \cdot\nabla \Psi_{\bm{\nu}} \right| \\
    &\leq  \left| a^* F_{\bm{\nu}} \right| + \left| \nabla a^* \cdot\nabla \Psi_{\bm{\nu}} \right| \\
    & \leq \left| a^* \right|  \left| F_{\bm{\nu}} \right| + \sum_{l=1}^d \left| \frac{\partial a^*}{\partial x_l}  \right| \left| \frac{\partial \Psi_{\bm{\nu}}}{\partial x_l}  \right| \\
    &\leq \left| a^* \right|  + \sum_{l=1}^d \left| \frac{\partial a^*}{\partial x_l}  \right|.
\end{align*}
Combining the above inequality with the conclusion \eqref{eq:eta*=0_infinite_norm} of \S\ref{s:eta*=0} (Step 2) yields
$$
    \left\|\Phi_{\bm{\nu}}    \right\|_{L^\infty}  \leq  \left|
    -\nabla\cdot\left(a \nabla_t \Psi_{\bm{\nu}} \right) \right|
    +  \left|
    -\nabla\cdot\left(a^* \nabla \Psi_{\bm{\nu}} \right) \right| 
 \leq  |e_{\bm{0}}| + \beta +\left\| a^* \right\|_{L^\infty} + \sum_{l=1}^d \left\| \frac{\partial a^*}{\partial x_l}  \right\|_{L^\infty}.
$$
This concludes the proof of Proposition~\ref{prop:compressible_eta}.

\subsection{Sparse recovery in bounded Riesz systems}

\label{app:RIP}

We start by recalling the definition of Restricted Isometry Property (RIP), which is a popular tool used to derive sufficient conditions for accurate and stable recovery in compressive sensing (see, e.g., \cite[Chapter 6]{foucart2013mathematical} and references therein).

\begin{definition}[Restricted Isometry Property (RIP)]
A matrix $M \in \mathbb{C}^{m \times N}$ has the Restricted Isometry Property (RIP) of order $k \in \mathbb{N}$ if there exists $0\leq\delta<1$ such that
\begin{equation}
    (1-\delta)\| \bm{z} \|_2^2 \leq \|A \bm{z}\|_2^2\leq (1+\delta)\|\bm{z}\|_2^2,
    \quad \forall \bm{z}\in \mathbb{C}^N, \; \|\bm{z}\|_0 \leq k.  
\end{equation}
The smallest $0\leq\delta<1$ such that the above condition holds is called the $k$-th restricted isometry constant $\delta_k$ of $M$. 
%where $\Sigma_s=\{\bm{z}\vert\ \|\bm{z}\|_0 \leq s\}$ represent all $s$-sparse vector in $\mathbb{R}^n$. %If a matrix $A \in \mathbb{C}^{m \times N}$ satisfies $\delta_s \leq \varepsilon$ for some $\varepsilon>0$ and a given sparsity $s$, then we say that $A$ satisfies $\mbox{RIP}(\varepsilon,s)$
\end{definition}
We now consider a sufficient lower bound on $m$ in order for the compressive Fourier collocation matrix $A$ to the RIP with high probability. This is an immediate consequence of the RIP theorem \cite[Theorem 2.3]{brugiapaglia2021sparse} for general random sampling in bounded Riesz systems and of Proposition~\ref{prop:compressible_eta} (that establishes that $\{\Phi_{\bm{\nu}}\}_{\bm{\nu} \in \Lambda}$ is a bounded Riesz system).
\begin{theorem}[RIP for High-dimensional Compressive Fourier Collocation] 
\label{thm:RIP}
Consider the same setting as in Proposition~\ref{prop:compressible_eta}. Then there exist a universal constant $c>0$ such that the following holds. Let  $\delta \in (1-b_\Phi/B_\Phi,1)$ and assume that the number of collocation points $m$ satisfies
$$
m \geq c \max\{B_\Phi^{-2},1\} K_\Phi^2 \xi^{-2} s L, 
$$
where
$$
L = \log^2(sK_\Phi^2\max\{B_\Phi^{-2},1\}\xi^{-2}) \log(eN) + \log(2\epsilon^{-1}),
$$
$\xi=\delta -1 +b_\Phi/B_\Phi>0$, and $N$ is the cardinality of the index set $\Lambda$. Then, the rescaled compressive Fourier collocation matrix $\frac{1}{\sqrt{B_{\Phi}}}A$, where $A$ is defined as in \eqref{eq:def_A_b}, has the $\mathrm{RIP}(\delta,s)$ with probability at least $1-\epsilon$.
%\begin{equation}
%\label{eq:high-probability}
%1-2\exp\left(\frac{-c_1 \min\{B_\Phi^2,1\}\xi^2m}{K_\Phi^2 s}\right).
%\end{equation}
%\purple{[To do: use probability of failure $\epsilon$ and longer log factor.]}

\end{theorem}

We now present a result showing that the RIP is sufficient to achieve accurate and stable recovery via OMP and QCBP. The following results combines \citet[Theorem 6.25]{foucart2013mathematical} (where the constant $1/13$ arises from taking into account column normalization -- see also \citet[Theorem 3.7]{brugiapaglia2021sparse}) and \citet[Theorem 2.1]{cai2013sparse} with $t=2$. For more details on the RIP-based analysis of OMP we also refer to the seminal paper \cite{zhang2011sparse} and \cite{cohen2017orthogonal}.
\begin{theorem}[RIP $\Rightarrow$ Accurate and stable recovery via OMP and QCBP]
\label{thm:RIP->recovery}
There exist universal constants $C_1,C_2 >0$ such that the following recovery guarantee holds for OMP and QCBP. Let $s\in \mathbb{N}$ and assume that $M \in \mathbb{C}^{m \times N}$ satisfies one of the following conditions:
\begin{equation}
\label{eq:RIP_condition}
\begin{cases}
\delta_{26s} < \frac{1}{13},  & \text{(OMP)}\\
\delta_{2s} < \frac{1}{\sqrt{2}}.  & \text{(QCBP)}
\end{cases}
\end{equation}
Then, for every $\bm{z} \in \mathbb{C}^N$ and every $\bm{e} \in \mathbb{C}^m$, the following holds. Let $\bm{b} = M \bm{z} + \bm{e}$ and $\hat{\bm{z}}$ be either the vector computed via $K=24s$ iterations of OMP (Algorithm \ref{alg:OMP}) or any solution to QCBP with $\gamma \geq \|\bm{e}\|_2$. Then, $\hat{\bm{z}}$ satisfies
$$
\| \bm{z} - \hat{\bm{z}} \|_2 
\leq C_1 \frac{\sigma_s(\bm{z})_1}{\sqrt{s}} + C_2 \|\bm{e}\|_2.
$$
\end{theorem}

\subsection{Proof of Theorem~\ref{thm:recovery_CFC}}
\label{app:recovery_CFC}

First, using Theorem~\ref{thm:RIP} we find a sufficient condition on $m$ that guarantees suitable RIP properties for the compressive Fourier collocation matrix $A$ needed to apply Theorem~\ref{thm:RIP->recovery}.

\paragraph{Step 1: Condition on $m$ $\Rightarrow$ RIP.} Let us focus on the OMP case first. Thanks to the sufficient condition \eqref{eq:suff_cond_on_Riesz_constants}, we have 
$$
1-\frac{b_{\Phi}}{B_{\Phi}} < \frac{0.98}{13}.
$$
Hence, since in order to apply Theorem~\ref{thm:RIP->recovery} we need the RIP with $\delta < 1/13$ and in Theorem~\ref{thm:RIP} we have the restriction $\delta > 1-b_{\Phi}/B_{\Phi}$, these two inequalities are compatible and the choice $\delta = 0.99/13$ is admissible. This leads to $$
\xi = \delta - 1 + \frac{b_{\Phi}}{B_{\Phi}} > \frac{0.01}{13}.
$$
Hence, the factor $\xi^{-2}$ in the sample complexity bound of Theorem~\ref{thm:RIP->recovery} can be dropped, up to choosing the universal constant $c$ large enough. We now justify the log factor in \eqref{eq:log_factor_CFC}. This is due to the fact that the cardinality bound \eqref{eq:cardinality_bound_Lambda} leads to
\begin{align*}
\log(\mathrm{e} N) 
& \leq \min\{\log(4\mathrm{e} n^5 16^d), \log(\mathrm{e}^3 n^{2+\log_2(d)})\}
\leq c_1 \min\{\log(n) + d, \log(2n) \log(2d)).
\end{align*}
This shows that the condition \eqref{eq:suff_cond_m_CFC} is sufficient to apply Theorem~\ref{thm:RIP->recovery} (with $s$ replaced by $26s$). The same argument can be applied to the case of QCBP up to modifying the constants accordingly. This shows that the matrix $A/\sqrt{B_{\Phi}}$ satisfies the RIP condition \eqref{eq:RIP_condition} with probability at least $1-\epsilon$.

\paragraph{Step 2: RIP $\Rightarrow$ Coefficient recovery.} We are now in a position to apply Theorem~\ref{thm:RIP->recovery}. Again, we discuss only the OMP case in detail. The proof for QCBP is similar. Recall the notation described in \S2.1. Assume the exact solution to the diffusion equation has the following expansions with respect to the spectral basis and to the $L^2$-orthonormal Fourier basis:
\begin{equation}
    u = \sum_{\bm{\nu} \in \mathbb{Z}^d} c_{\bm{\nu}} \Psi_{\bm{\nu}}
    = \sum_{\bm{\nu} \in \mathbb{Z}^d} \tilde{c}_{\bm{\nu}} F_{\bm{\nu}}.
\end{equation}
Here the $\bm{c}, \bm{\tilde{c}}\in \ell^2(\mathbb{Z}^d)$ are the coefficients of the spectral basis (with rescaling) and orthogonal Fourier basis respectively. They satisfy the relation
$$
c_{\bm{\nu}} = (2\pi)^2 \|\bm{\nu}\|_2^2 \,\tilde{c}_{\bm{\nu}}, \quad \forall \bm{\nu} \in \mathbb{Z}^d.
$$
Moreover, their truncated versions satisfy
\begin{equation}
    \bm{c}_\Lambda = D \bm{\tilde{c}}_{\Lambda},
\end{equation}
where $D=(2\pi)^2 \diag\{{\| \bm{\nu}_1 \|}_2^2, \dots, {\| \bm{\nu}_N\|}_2^2 \}$. Recall that the truncated solution on a finite multi-index set $\Lambda$ is denoted by $u_\Lambda$. Then 
$$
u_\Lambda = \sum_{\bm{\nu} \in \Lambda} c_{\bm{\nu}} \Psi_{\bm{\nu}}
= \sum_{\bm{\nu} \in \Lambda} \tilde{c}_{\bm{\nu}} F_{\bm{\nu}}.
$$
We recall from \eqref{eq:linear_system_c_Lambda} that 
$$
A \bm{c}_{\Lambda} = \bm{b} + \bm{e}.
$$
Therefore, thanks to Theorem~\ref{thm:RIP->recovery}, we have
$$
\|\bm{c}_\Lambda - \hat{\bm{c}}\|_2 \leq
C_1\frac{\sigma_s(\bm{c}_\Lambda)_1}{\sqrt{s}}+ C_2 \|\bm{e}\|_2.
$$

\paragraph{Step 3: Coefficient recovery $\Rightarrow$ Solution recovery.} To obtain \eqref{eq:CFC_L2_Laplacian_bound}, we first apply the triangle inequality
\begin{align*}
\|\Delta(u-\hat{u}) \|_{L^2} 
& \leq \|\Delta (u-u_\Lambda)\|_{L^2} + \|\Delta (\hat{u}-u_\Lambda)\|_{L^2}.
\end{align*}
Then, recalling that $-\Delta \Psi_{\bm{\nu}} = F_{\bm{\nu}}$, we note that
$$
\|\Delta (u_\Lambda- \hat{u})\|_{L^2} = 
\|\bm{c}_\Lambda - \hat{\bm{c}}\|_2 \leq \frac{C_1 \sigma_s(\bm{c}_\Lambda)_1}{\sqrt{s}} + C_2 \| \bm{e}\|_2.
$$
Finally, recall from \eqref{eq:error_norm_bound} that 
$\|\bm{e}\|_2 \leq \|\nabla \cdot (a \nabla (u - u_{\Lambda}))\|_{L^\infty}$. Combining the above inequalities yields \eqref{eq:CFC_L2_Laplacian_bound}.

Let us now show \eqref{eq:CFC_L2_bound}. Using the triangle inequality, we see that
$$
\|u-\hat{u} \|_{L^2} \leq \|u - u_\Lambda\|_{L^2} + \|u_\Lambda -  \hat{u} \|_{L^2}
$$
%The first term on the right is the truncation error. 
%$$
% \| u - u_\Lambda \|_{L^2} = \| \tilde{\bm{c}} - \tilde{\bm{c}}_\Lambda \|_2
%    = \sqrt{\sum_{\nu \not\in \Lambda} \tilde{c}_{\nu}^2}
%$$
Using that $D_{ii}\geq 4\pi^2$ (because $\bm{0} \notin\Lambda$), the second term can be bounded as
$$
\| u_\Lambda- \hat{u} \|_{L^2} = \|D^{-1}(\bm{c}_\Lambda - \hat{\bm{c}})\|_2
\leq \frac{1}{4\pi^2} \left( \frac{C_1 \sigma_s(\bm{c}_\Lambda)_1}{\sqrt{s}} + C_2 \| \bm{e}\|_2
\right)
$$
This, combined with \eqref{eq:error_norm_bound}, implies the $L^2$-error bound \eqref{eq:CFC_L2_bound} and concludes the proof.

\subsection{Further details on numerical experiments}
\label{a:details_numerics}

\paragraph{Hardware and software specifications.} All the experiments have been performed on a Dell G3 3590 with 8GB of RAM and 2.40GHz Intel i5-9300H processor, using MATLAB\textregistered\ 2019b.

\paragraph{Computation of \eqref{eq:def_A_b}.} The right-hand side $\bm{b}$ of the linear system \eqref{eq:def_A_b} is evaluated using the exact solution via sixth-order finite difference scheme. The entries of the matrix $A$ are instead computed explicitly. 

\paragraph{Recovery via OMP, and QCBP.} For OMP (Algorithm~\ref{alg:OMP}), we use $K=m/2$ iterations. The QCBP optimization program \eqref{eq:QCBP} is numerically solved using the CVX MATLAB\textregistered\ package \cite{grant2008graph,grant2014cvx} with the MOSEK solver and the \texttt{cvx\_precision} parameter always set as \texttt{'default'}.  The tuning parameter $\eta$ of QCBP is chosen according to an ``oracle'' approach. Namely, a reference solution $\tilde{\bm{c}} \approx \bm{c}_\Lambda$ is pre-computed using a large number of collocation points ($m=4|\Lambda|$). Then, the tuning parameter $\eta$ is set as $\eta = \|A\tilde{\bm{c}} - \bm{b}\|_2$. We note that this idealized strategy is not recommended in practice since it requires a large number of collocation points. In practice, one could tune $\eta$ via cross-validation or, even better, replace QCBP with a ``noise blind'' decoder such as the square-root LASSO,  which does not require any \emph{a priori} estimate of the unknown error norm $\|\bm{e}\|_2$ (see \cite{adcock2019correcting}). 

\paragraph{Randomization of experiments and visualization.} Due to the randomized nature of compressive Fourier collocation, we consider 25 random runs for each test.  Moreover, the number of sampling points $m$ is always a power of 2 in our experiments. The curves represent the sample geometric mean of the errors, and the size of the lightly shaded areas is the corrected standard geometric mean. For more details on the visualization strategy, we refer to \cite[Appendix A.1.3]{adcock2022sparse}

%USE THE BELOW OPTIONS IN CASE YOU NEED AUTHOR YEAR FORMAT.
\bibliographystyle{abbrvnat}
\bibliography{biblio}

\end{document}